\documentclass[12pt]{amsart}
\newcommand{\et}{\quad\mbox{and}\quad}

\newcommand{\bC}{\mathbb{C}}

\newcommand{\bQ}{\mathbb{Q}}
\newcommand{\bR}{\mathbb{R}}
\newcommand{\bZ}{\mathbb{Z}}
\newcommand{\bF}{\mathbb{F}}
\newcommand{\bN}{\mathbb{N}}

\newcommand{\bI}{\mathbb{I}}

\newcommand{\GL}{\mathrm{GL}}
\newcommand{\M}{\mathrm{M}}

\newcommand{\Kbar}{{\bar{K}}}

\newcommand{\disp}{\displaystyle}

\newcommand{\Adj}{\mathrm{Adj}}
\def\disp{\displaystyle}
 \numberwithin{equation}{section}
\begin{document}
\baselineskip=17pt
\title[Equations diff\'erentielles $p$-adiques et S\'eries Gevrey  arithm\'etiques]
{Equations diff\'erentielles $p$-adiques et S\'eries Gevrey
arithm\'etiques}
\author{Said Manjra et Salah-eddine Remmal}
\address{
   D\'epartement de Math\'ematiques\\
   Universit\'e d'Ottawa\\
   585 King Edward\\
   Ottawa, Ontario K1N 6N5, Canada}
\email{manjra@math.net}
\address{
   D\'epartement de Math\'ematiques\\
   Facult\'e des sciences de F\`es\\
   BP 1796. Atlas-F\`es. \\
   F\`es. Maroc}
\email{Remmal@math.net} 

\thanks{Travail partiellement subventionn\'e par le CRSNG}
\maketitle

\section{ Introduction}

Cet article est consacr\'e \`a une \'etude $p$-adique de
 l'\'equation diff\'erentielle alg\'ebrique de degr\'e minimal annulant une s\'erie
Gevrey arithm\'etique d'ordre $-1$ (c'est-\`a-dire, une
$E$-fonction). Plus pr\'ecis\'ement, on y
d\'emontre une conjecture d'Yves Andr\'e (conjecture 4.7 de [A2]).\\

Fixons d'abord quelques notations. Soit $K$ un  corps de nombres et
soit $V_{0}$ l'ensemble de toutes les places finies $v$ de $K$. Pour
chaque $v\in V_{0}$ au-dessus d'un nombre premier $p=p(v)$, on
normalise la valeur absolue $v$-adique de fa\c{c}on que
$|p|_{v}=p^{-1}$ et on note $\pi_{v}=p^{\frac{-1}{p-1}}$. On
d\'esigne aussi par $K_{v}$ le compl\'et\'e de $K$ pour cette valeur
absolue, et on fixe un plongement $K\hookrightarrow {\bC}$. Pour
tout r\'eel $r>0$, le \emph{rayon de convergence g\'en\'erique}
$R_{v}(\phi,r)$ d'un op\'erateur diff\'erentiel $\disp\phi \in
K[x,d/dx]$ est, par d\'efinition, le rayon de convergence, limit\'e
sup\'erieurement par $r$, d'une base de solutions de $\phi$ au
voisinage d'un point g\'en\'erique $v$-adique de valeur absolue $r$.\\

Une s\'erie formelle $\sum_{n\ge 0}a_nx^n$ est dite \emph{s\'erie
Gevrey arithm\'etique} d'ordre $s\in {\bQ}$ si ses coefficients
$a_n$ sont des nombres alg\'ebriques et s'il existe une constante
$C>0$ telle que les conjugu\'es du nombre alg\'ebrique $a_n/(n!)^s$
soient de valeur absolue inf\'erieure \`a $C^n$, et que, pour tout
$n\in {\bN}$, le d\'enominateur commun $d_n$ des nombres
$a_0=a_0/(0!)^s, \ldots, a_n/(n!)^s$ soit inf\'erieur \`a $C^n$.
Cette classe englobe les s\'eries hyperg\'eom\'etriques
g\'en\'eralis\'ees confluentes ou non, \`a param\`etres rationnels.
L'ensemble de telles s\'eries est not\'e $\overline{\bQ}\{x\}_s$.
Une s\'erie formelle \`a coefficients alg\'ebriques est dite
holonome si elle est solution d'une \'equation diff\'erentielle
lin\'eaire \`a coefficients dans $\overline{\bQ}(x)$. Les
$G$-fonctions (resp. les $E$-fonctions) sont par d\'efinition les
\'el\'ements holonomes de
$\overline{\bQ}\{x\}_0$ (resp. de $\overline{\bQ}\{x\}_{-1}$). \\

La th\'eorie des $G$-fonctions  a connu un essor remarquable au
cours des derni\'eres d\'ecennies, principalement gr\^ace aux
travaux de Chudnovsky, Bombieri et Andr\'e. Leurs contributions ont
mis en \'evidence certaines connexions avec la transcendance et la
g\'eom\'etrie arithm\'etique. En particulier, Chudnovsky a prouv\'e
que l'op\'erateur diff\'erentiel  $\phi$ de $K[x,d/dx]$ d'ordre
minimal (en $d/dx$) annulant une $G$-fonction satisfait la condition
dite de Galochkin (cf. [CC], [A1, VI]). Un op\'erateur
diff\'erentiel de $K[x,d/dx]$ qui remplit cette condition est
appel\'e un \emph{$G$-op\'erateur}.  Andr\'e, plus tard, a \'etabli
l'\'equivalence entre cette condition et celle de Bombieri, \`a
savoir qu'un op\'erateur diff\'erentiel $\phi$ de $K[x,d/dx]$ est un
$G$-op\'erateur si et seulement s'il satisfait $\prod_{v\in
V_0}R_{v}(\phi,1)\ne 0$ [A1, IV 5.2]. Gr\^ace au th\'eor\`eme de
Katz [Ka, 13.0], on v\'erifie que les $G$-op\'erateurs n'ont que des
singularit\'es r\'eguli\`eres \`a exposants rationnels. De plus,
toutes les s\'eries formelles qui interviennent dans leurs
solutions, en tout point fini ou
infini, sont des $G$-fonctions [A1, V].\\

L'\'etude des s\'eries Gevrey arithm\'etiques d'ordre non nul
holonomes a \'et\'e r\'ecemment d\'evelopp\'ee par Y. Andr\'e dans
[A2]. En particulier, ce dernier a introduit une nouvelle classe
d'op\'erateurs diff\'erentiels appel\'es $E$-op\'erateurs. Par
d\'efinition, un $E$-op\'erateur est la transform\'ee de
Fourier-Laplace d'un $G$-op\'erateur. Cette appellation est
justifi\'ee par le fait que l'op\'erateur diff\'erentiel de degr\'e
minimal en $x$ de $K[x,d/dx]$ annulant une $E$-fonction est la
transform\'ee de Fourier-Laplace d'un $G$-op\'erateur [A2, 4.2].
\'Egalement, les s\'eries formelles qui interviennent dans les
solutions  en $0$ d'un $E$-op\'erateur sont toutes des $E$-fonctions
[A2, 4.3]. L'importance de ces op\'erateurs vient du fait que, si
$y\in \overline{\bQ}\{x\}_s$ pour un certain rationnel $s$ non nul,
alors $y(x^{-s})$ est solution d'un $E$-op\'erateur $\psi$ [A2,
6.1].  Un tel r\'esultat, conbin\'e aux propri\'et\'es de ces
op\'erateurs, fournit une bonne connaissance des propri\'et\'es de
l'op\'erateur d'ordre minimal annulant un \'el\'ement de
$\overline{\bQ}\{x\}_{s\ne 0}$ (Th\'eor\`eme de dualit\'e [A2]). Ces
r\'esultats  permettent, en particulier, de retrouver certains
th\'eor\`emes fondamentaux de la th\'eorie de transcendance:
\emph{th\'eor\`eme de Lindemann-Weierstrass et th\'eor\`eme de
Siegel-Shidlovskii} [A3]. Andr\'e a aussi conjectur\'e une
caract\'erisation $p$-adique des $E$-op\'erateurs analogue \`a celle
de Bombieri pour les $G$-op\'erateurs [A2, 4.7]. En r\'eponse \`a
cette
conjecture, nous d\'emontrons ici le r\'esultat suivant: \\

\ \\{\bf Th\'eor\`eme 1.1.}\emph{ Soit $\psi\in K[x,d/dx]$. Alors
$\psi$ est un $E$-op\'erateur si et seulement si  les conditions
suivantes sont satisfaites:
\begin{itemize}
 \item[(1)] les pentes du polygone de Newton $($au sens de Ramis$)$
de $\psi$ sont dans $\{-1,0\}$,
 \item[(2)] ${\prod_{v\in V_{0}} R_{v}(\psi,\pi_{v})\pi_{v}^{-1}\ne 0}$.
\end{itemize}
}\

Notons que la conjecture pr\'ecise d'Andr\'e remplace la condition
(1) par la condition plus faible que $\phi$ n'ait de singularit\'es
qu'en $0$ et en l'infini.  Cependant, on verra par un exemple que
ceci ne suffit pas et qu'il y faut la condition plus forte que nous
donnons sur le polygone de Newton-Ramis (voir l'exemple 4.6). \\

La d\'emonstration de ce th\'eor\`eme est bas\'ee essentiellement
sur les propri\'et\'es locales des $G$-op\'erateurs (Corollaire
3.2.1) par le biais d'une transform\'ee de Laplace formelle
(nouvelle) et  d'une  g\'en\'eralisation du principe de transfert de
Christol [MR]  combin\'es \`a quelques propri\'et\'es du
polygone de Newton-Ramis.\\

Cet article est organis\'e ainsi:\\

La section 2 est consacr\'ee aux notations et r\'esultats
pr\'eliminaires. En particulier, on y donne un crit\`ere
arithm\'etique local pour la rationalit\'e des exposants (Corollaire
2.3.3), et on y \'etablit l'\'equivalence entre conditions
``Bombieri'' et ``local Bombieri'' dans un cadre plus g\'en\'eral
(Th\'eor\`eme 2.4.2). On d\'efinit ensuite le polygone de Newton au
sens de Ramis et on cite certaines de ses propri\'et\'es (\S2.5). On
d\'efinit aussi une transform\'ee de Fourier-Laplace
g\'en\'eralis\'ee qui \'etend \`a la fois la transform\'ee de
Fourier-Laplace usuelle et celle de Dwork normalis\'ee (\S2.6). Dans
la section 3, on rappelle quelques propri\'et\'es des
$G$-op\'erateurs et $E$-op\'erateurs. Puis, dans la section 4, on
d\'emontre une premi\`ere moiti\'e du th\'eor\`eme 1.1, \`a savoir
que les conditions (1) et (2) de ce th\'eor\`eme sont n\'ecessaires.
On introduit, dans la section 5, une transformation de Laplace
formelle g\'en\'eralis\'ee $\mathcal{L}_\Lambda^\tau$, par rapport
\`a une matrice $\Lambda$ carr\'ee $n\times n$ \`a coefficients dans
un corps de nombres $K$ et \`a un nombre alg\'ebrique $\tau\in K$,
agissant sur les matrices de la forme de $A(x)x^\Lambda$, o\`u
$A(x)$ est une matrice $m\times n$ \`a coefficients dans
$K[[x,1/x]]$. Cette transformation satisfait des propri\'et\'es de
commutation avec la d\'erivation $d/dx$ et avec la transform\'ee de
Fourier-Laplace g\'en\'eralis\'ee qui sont analogues au cas
classique (Proposition 5.1).  De plus on d\'emontre des relations de
comparaison entre les rayons de convergence des coefficients de
$A(x)x^\Lambda$ et ceux de son image $\mathcal{L}_\Lambda^\tau(
A(x)x^\Lambda )$ (Proposition 5.2.2).  On combine tous ces outils
dans la section 6 pour compl\'eter la d\'emonstration du
th\'eor\`eme 1.1.\\
\ \\{\bf Remerciements.} \emph{Nous tenons \`a remercier G. Christol
pour son soutien, D. Roy pour ses suggestions dans la pr\'esentation
du \S5 et   Y. Andr\'e pour ses remarques.}
\section{ Notations et r\'esultats \'el\'ementaires}
\ \\
{\bf 2.1 Modules diff\'erentiels}\\

Soit $\mathcal K$ un corps commutatif muni d'une d\'erivation
$\partial$, soit $K$ le corps des constantes de $\partial$ dans
$\mathcal K$ et soit $n\ge 1$ un entier.
  Un $\mathcal K$-module
diff\'erentiel $\mathcal M$ est un module libre de rang $n$ sur
$\mathcal K$ muni
  d'un  $K$-endomorphisme $\nabla$
   de $\mathcal M$ qui v\'erifie pour tout $m\in\mathcal M$ et tout
$a\in\mathcal
   K$ la condition $\nabla(am)=a\nabla(m)+\partial(a)m.$ A chaque base
$\{e_{i}\}$ de $\mathcal
   M$ sur $\mathcal K$
  correspond une matrice $A=(A_{ij})\in {\M}_n({\mathcal K})$ v\'erifiant
  $$\nabla(e_{i})=\sum_{j=1}^{n}A_{ij}e_{j}$$ appel\'ee
  repr\'esentation de la d\'erivation $\partial$ dans la base $\{e_{i}\}$
  (ou simplement matrice associ\'ee \`a $\mathcal M$) ainsi qu'un syst\`eme
   diff\'erentiel  $\partial X=AX$ o\`u $X$ d\'esigne un vecteur colonne
   $n\times 1$ ou encore une matrice $n\times n$.
    Un changement de base dans $\mathcal M$ se
   traduit par la donn\'ee d'une matrice $Y\in {\GL}_n({\mathcal K})$ telle
que
   $Y[A]:=YAY^{-1}+\partial(Y)Y^{-1}$ soit la repr\'esentation de $\partial$
dans
   la nouvelle base. La donn\'ee  de $\mathcal M$ correspond aussi \`a celle
   d'une \'equation diff\'erentielle
    lin\'eaire via le lemme du vecteur cyclique (cf. $[$Ma, III$]$).\\

Invers\'ement, si on dispose d'un op\'erateur diff\'erentiel $\phi=
\disp\sum_{i=0}^{n} a_{i}\partial^{i}\in {\mathcal K}[\partial]$
avec $ a_{n}\ne 0$, on peut lui associer le ${\mathcal K}$-module
diff\'erentiel
    ${\mathcal M}_{\phi}={\mathcal K}[\partial]/{\mathcal K}[\partial]\phi$
de rang $n$ qui correspond au syst\`eme $$
\partial X=A_{\phi}X\;\;\mbox{o\`u}\;\; A_{\phi}:=
\left (
\begin{array}{c}
0\\0\\\vdots\\0\\-\frac{a_0}{a_{n}}
\end{array}
\begin{array}{c}
1\\0\\\vdots\\0\\\frac{a_1}{a_{n}}
\end{array}
\begin{array}{c}
0\\1\\\vdots\\0\\-\frac{a_2}{a_{n}}
\end{array}
\begin{array}{c}
\ldots\\\ldots\\\ddots\\ \\\ldots
\end{array}
\begin{array}{c}
0\\0\\\vdots\\1\\-\frac{a_{n-1}}{a_{n}}
\end{array}
\right).
$$

On associe  aussi \`a $\phi$ l'op\'erateur adjoint $\phi^{*}=
\disp\sum_{i=0}^{n}(-\partial)^{i}a_{i}$. On v\'erifie que
$-^{t}A_{\phi}$ est associ\'ee \`a ${\mathcal
M}_{\phi^{*}}={\mathcal K}[\partial]/{\mathcal
K}[\partial]\phi^{*}$. Plus g\'en\'eralement, $A$ est une matrice
associ\'ee \`a ${\mathcal M}_{\phi}={\mathcal K}[\partial]/{\mathcal
K}[\partial]\phi$,
     si et seulement si  $-^{t}A$ est associ\'ee \`a  ${\mathcal
     M}_{\phi^{*}}$. Ceci r\'esulte du fait que pour tout $Y\in
     {\GL}_n({\mathcal K})$, on a
     \begin{eqnarray*}
        -^{t}(Y[A_{\phi}]) &=&
^{t}Y^{-1}(-^{t}A_{\phi})^{t}Y-^{t}Y^{-1}\partial(^{t}Y)
=^{t}Y^{-1}(-^{t}A_{\phi})^{t}Y+\partial(^{t}Y^{-1})^{t}Y\\
                                     &=&
                                     (^{t}Y^{-1})[-^{t}A_{\phi}].
\end{eqnarray*}

Soient  $A_1,\ldots,A_\ell$  des  matrices carr\'ees \`a
coefficients dans $\mathcal{K}$.
  On note par $\mathcal{M}_{A_i}$ le module diff\'erentiel associ\'e \`a  la
matrice $A_i$
  pour $i=1,\ldots,\ell$, et par
  $\oplus_{1\le i\le \ell}A_{i}$ la matrice diagonale par blocs
$$\oplus_{1\le i\le \ell}A_{i}= \left (
\begin{array}{c}
A_1\\ \ \\ \ \\ \
\end{array}
\begin{array}{c}
\ \\A_2 \\ \ \\ \
\end{array}
\begin{array}{c}
\ \\ \ \\ \ddots\\ \
\end{array}
\begin{array}{c}
\ \\ \ \\ \ \\ A_\ell
\end{array}
\right)$$ avec les blocs $A_1,\ldots,A_\ell$ sur la diagonale. De
fa\c{c}on similaire, on d\'efinit $\oplus_{1\le i\le
\ell}\mathcal{M}_{A_{i}}:=\mathcal{M}_{\oplus_{1\le i\le
\ell}A_{i}}$. Si $A_1,\ldots,A_\ell\in {\M}_n(\mathcal{K})$, on pose
$\otimes_{1\le i\le
\ell}\mathcal{M}_{A_{i}}:=\mathcal{M}_{A_1+\ldots +A_\ell}$. On
notera par ${\bI}_n$ la matrice identit\'e de
${\M}_n(\mathcal{K})$.\\ \
\\
{\bf 2.2 Rayon de convergence au voisinage d'une singularit\'e}\\

On consid\`ere le corps diff\'erentiel ${\mathcal K}=K(x)$ muni de
la d\'erivation $\partial=d/dx$. Soit $\mathcal M$ le $K(x)$-module
diff\'erentiel de rang $n$ associ\'e \`a une matrice $A(x)$ \`a
coefficients dans $K(x)$. On dit que $\mathcal M$ est r\'egulier en
$0$, s'il existe une extension finie $K'$ de $K$, une matrice
$Y_A(x)\in {\GL}_n(K'((x)))$, appel\'ee \emph{matrice de
r\'eduction} de $A(x)$ (ou de $\mathcal M$)
  en $0$ et une matrice $\Gamma\in
  {\M}_n(K')$ telles que $Y_A(x)[A(x)]=\Gamma/x$. Dans ce cas,
  $Y_A(x)^{-1}x^{\Gamma}$ est solution du syst\`eme
    $dX/dx=A(x)X$. On v\'erifie que les valeurs propres de $\Gamma$
modulo $\mathbb{Z}$ ne d\'ependent que de $\mathcal M$. On les
appelle les \emph{exposants} de  $\mathcal M$ en $0$. Quitte \`a
faire un changement de bases dans $\mathcal{M}$, on peut supposer
que la matrice $\Gamma$ est sous forme de Jordan.\\

En particulier, si $A(x)=A_{\phi}$ pour un op\'erateur
diff\'erentiel $\phi$ de $K(x)[\partial]$, la matrice $Y_A(x)$
sera not\'ee $Y_\phi$ et appel\'ee matrice de r\'eduction de $\phi$.\\

On dit que ${\mathcal M}$ est r\'egulier en  l'infini, si le
$K(x)$-module ${\mathcal M}^{\infty}$
     obtenu \`a partir de $\mathcal M$ par le changement de variable
$x\longrightarrow 1/x$ est r\'egulier  en $0$.
     Cela signifie qu'il existe une extension
finie $K'$ de $K$, une matrice $Y_{\infty}(x)\in {\GL}_n(K'((x)))$,
appel\'ee matrice de r\'eduction de $A$ (ou de $\mathcal M$)  en
l'infini et une matrice $\Gamma_{\infty}\in
  {\M}_n(K')$   telles que
   $Y_{\infty}(x^{-1})[A(x)]=-\Gamma_{\infty}/x$. Les exposants de
${\mathcal M}$ en l'infini sont les valeurs propres de
$\Gamma_{\infty}$ modulo $\mathbb{Z}$.
   Dans ce cas, $Y_{\infty}(x^{-1})(x^{-1})^{\Gamma_{\infty}}$ est solution
du syst\`eme  $dX/dx=A(x)X$.\\

Soit $\phi\in K(x)[\partial]$. Par extension,
  on attribuera \`a $\phi$ les propri\'et\'es que
${\mathcal M}_{\phi}$ poss\`ede. On constate alors, depuis \S 2.1,
que $\phi$ est r\'egulier en $0$ (resp. l'infini) si et seulement si
$\phi^{*}$ est  r\'egulier au  m\^eme point, auquel cas, les
exposants de $\phi^{*}$ en $0$ (resp. en l'infini) sont les m\^emes
que ceux de $\phi$  mais de signes
oppos\'es.\\

On d\'esignera par $\overline{\phi}$ l'op\'erateur diff\'erentiel
obtenu \`a partir de $\phi$ par le changement de variable
$x\longrightarrow -x$ et par ${\mathcal M}_{\alpha}$ le
$K(x)$-module  diff\'erentiel (r\'egulier en $0$) associ\'e \`a
l'\'equation $xd/dx-\alpha$ o\`u $\alpha$ est un \'el\'ement de
  $K$.\\

Dans la suite, on pourra toujours supposer $K$ suffisamment grand de
sorte qu'on puisse prendre $K'=K$. On fixera  pour toute place finie
$v$ de $V_{0}$, un plongement $K\hookrightarrow
K_{v}\hookrightarrow{\bC}_{p(v)}\hookrightarrow\Omega_{p(v)}$, o\`u
$\Omega_{p(v)}$ est une extension de ${\bC}_{p(v)}$ valu\'ee
compl\`ete  et alg\'ebriquement close, dont le groupe des valeurs
absolues est ${\bR}_{\ge 0}$ et dont le corps r\'esiduel est une
extension stricte du corps r\'esiduel ${\bF}_{p(v)^{\infty}}$ de
${\bC}_{p(v)}$. Pour tout r\'eel $r>0$, on notera par $t_{v,r}$ un
point g\'en\'erique dans $\Omega_{p(v)}$ de valeur absolue $r$. On
notera par $S$ l'ensemble de tout les nombres premiers de
${\bZ}_{>0}$. On dira qu'une propri\'et\'e est v\'erifi\'ee pour
presque tout $v$, si elle est v\'erifi\'ee pour tout $v\in V_0$
except\'e
peut-\^etre pour un nombre fini de $v$.\\

Soient  maintenant $v$ une place finie fix\'ee de $V_{0}$ et $I$ un
intervalle de ${\bR}_{>0}$. On notera par $\mathcal{A}_{v}(I)$
l'anneau des fonctions analytiques dans la couronne
$\mathcal{C}_v(I):=\{a\in {\bC}_{v(p)}\;|\; |a|_v\in I\}$ :
$$\mathcal{A}_{v}(I)=\{\sum_{i\in{\bZ} }a_ix^i\in
{\bC}_{p(v)}[[x,1/x]]\;|\;\lim_{i\to\mp\infty}|a_i|_v r^i=0,
\forall\; r\in I\},$$ et par $\mathcal{H}_{v}(I)$ le compl\'et\'e de
l'anneau des fractions  rationelles $f$ de ${\bC}_{p(v)}(x)$ n'ayant
pas de p\^ole dans $\mathcal{C}_v(I)$ pour la norme
$\|f\|_{v,I}:=\sup_{r\in I}|f(t_{v,r})|_{v}$. On sait que
$\mathcal{H}_{v}(I)\subseteq \mathcal{A}_{v}(I)$, avec \'egalit\'e
si $I$ est ferm\'e. \\

Si $n\in {\bZ}_{\ge 0}$ et $Y\in {\M}_n(K((x)))$, on pose, s'il
existe, $r_{v}(Y):=\sup\{r> 0\; |\; Y\in
{\M}_n(\mathcal{A}_{v}(]0,r[))\}$, et
$R_{v}(Y):=\min(r_{v}(Y),r_{v}(Y^{-1}))$  si en plus $Y\in
{\GL}_n(K((x)))$. Cette  derni\`ere notation est justif\'ee par le
r\'esultat suivant de F. Baldassarri [Ba, III].
\\ \ \\{\bf Proposition 2.2.1.} \emph{Si $Y(x)$ est une
matrice de r\'eduction d'un $K(x)$-module $\mathcal M$ en $0$ ou en
l'infini, alors
$R_{v}(Y)$  est non nul.}\\

Le symbole de Pochhammer $(\alpha)_i$ est d\'efini pour tout
$\alpha\in K$ et tout entier $i\in {\bZ}_{\ge 0}$ par $(\alpha)_0=1$
et $(\alpha)_i=\alpha(\alpha+1)\ldots(\alpha+i-1)$ pour $i>0$ . Avec
cette notation, le th\'eor\`eme 2 de \cite{Cl} donne
aussit\^ot:\\ \ \\
\textbf{Lemme 2.2.2.} \emph{Soit $v$ une place finie de $K$. Pour
tout entier $\alpha\ge 1$ ou tout nombre rationel $\alpha$ non
entier de d\'enominateur premier \`a $p(v)$, on a} $ \lim_{i\to
\infty} |(\alpha)_i|_v^{1/i} = \pi_v.$ \emph{En particulier, si}
$\alpha=1$, \emph{on a}
$\disp\lim_{i\to\infty} |i!|_v^{1/i}=\pi_v$.\\ \ \\
{\bf 2.3 Rayon de convergence g\'en\'erique}\\

Soient  $n$ un entier positif, $v$ une place finie de $V_0$, $I$ un
intervalle de ${\bR}_{>0}$ et $A$ une  matrice  $n\times n$ \`a
coefficients dans $\mathcal{H}_{v}(I)$. Pour tout r\'eel $r$ de
$\overline{I}$ (l'adh\'erence de $I$ dans ${\bR}_{>0}$), le  rayon
de convergence g\'en\'erique $R_v(\mathcal{M}_A,r)$ du module
diff\'erentiel $\mathcal{M}_A$ associ\'e \`a $A$ est, par
d\'efinition, le rayon de convergence, limit\'e sup\'erieurement par
$r$, d'une solution $\mathcal{U}_{t_{v,r}}$ de
${\GL}_n(\Omega_{p(v)}[[x-t_{v,r}]])$ du syst\`eme $dU/dx=AU$ au
voisinage du point g\'en\'erique $t_{v,r}$. D'apr\`es [CR, 9.2.7],
si $Y\in {\GL}_n(\mathcal{H}_{v}(I)))$, on a pour tout $r$ de $I$
\begin{equation}
\label{chg} R_v(\mathcal{M}_{Y[A]},r)=R_v(\mathcal{M}_A,r).
\end{equation}

Du fait qu'on a pour tout sous-intervalle ferm\'e non vide $J$ de
$I$, $\mathcal{H}_{v}(J)=\mathcal{A}_{v}(J)$ , la d\'efinition  du
rayon de convergence g\'en\'erique ci-dessus et l'\'egalit\'e
\eqref{chg} s'\'etendent, pour tout $r\in I$, au cas o\`u $A\in
{\M}_n(\mathcal{A}_{v}(I))$,
$Y\in{\GL}_n(\mathcal{A}_{v}(I))$.\\

Supposons dor\'enavant que  $A\in {\M}_n(K(x))$. On d\'eduit du
th\'eor\`eme 2 de [Po] que le graphe de la fonction $\rho\mapsto
\log(R_v(\mathcal{M}_A,\exp(\rho)))$ est un polygone,  ayant un
nombre fini de c\^ot\'es \`a pentes rationnelles, appel\'e
\emph{polygone de convergence} de
  $\mathcal{M}_A$, sur ${\bR}$ que l'on note  par ${\mathcal
P}_{v}$. En particulier, si $A\in {\M}_n(\mathcal{H}_{v}(I))$, alors
la partie  de ${\mathcal P}_{v}$ au-dessus de l'intervalle
$\log(\overline{I}):=\{\log(\rho)\;|\; \rho\in \overline{I}\}$ est
concave. De plus, si $A=A_{\phi}$ pour un certain op\'erateur
diff\'erentiel $\phi$ de $K(x)[d/dx]$, on v\'erifie que les deux
fonctions $R_v(\mathcal{M}_\phi,.)$ et $R_v(\phi,.)$ co\"{\i}ncident
sur ${\bR}_{>0}$. Le polygone ${\mathcal P}_{v}$ de
$\mathcal{M}_\phi$ est alors dit aussi polygone de convergence de
$\phi$.
\\ \ \\{\bf Th\'eor\`eme 2.3.1} (de finitude) $[$MR,
3.1$]$. \emph{L'ensemble des pentes des
  polygones de convergence d'un op\'erateur diff\'erentiel $\phi$ de
$K(x)[d/dx]$, associ\'es
aux  diff\'erentes places finies de $K$, est fini.}\\

Ce th\'eor\`eme jouera un r\^ole essentiel dans la sous-section
suivante et par cons\'equent dans la d\'emonstartion du r\'esultat
principal du pr\'esent
article.\\

Supposons toujours que $A\in {\M}_n(K(x))$ et $v\in V_0$. Si $r$ est
un r\'eel $>0$ et $\mathcal{U}_{t_{v,r}}$ est une solution du
syst\`eme $dU/dx=AU$ au voisinage de $t_{v,r}$, alors le  rayon de
convergence de $\mathcal{U}_{t_{v,r}}^{-1}$ est sup\'erieur ou
\'egale \`a $R_v(\mathcal{M}_A,\rho)$ (cf. [CR, 9.2.5]). En
combinant cette  remarque avec le fait que, si $B$ est une fraction
rationnelle de $K(x)$ et  $\mathcal{V}_{t_{v,r}}$ est une solution
du syst\`eme $dV/dx=BV$ au voisinage de $t_{v,r}$, alors le produit
$\mathcal{U}_{t_{v,r}}\mathcal{V}_{t_{v,r}}$ est une solution du
syst\`eme  $dW/dx=(B{\bI}_n+A)W$ au voisinage $t_{v,r}$, on trouve
que pour tout $r>0$,
\begin{equation}
\label{tns} R_v(\mathcal{M}_{A}\otimes
\mathcal{M}_{B},r)\ge\min\{R_v(\mathcal{M}_{A},r),R_v(\mathcal{M}_{B},r)\},\;\;\text{avec
\'egalit\'e si }\;\; R_v(\mathcal{M}_{A},r)\ne
R_v(\mathcal{M}_{B},r).
\end{equation}
Par ailleurs, il est facile de voir que, si $B$ est une matrice
carr\'ee \`a coefficients dans $K(x)$, on a pour tout $r>0$,
\begin{equation}
\label{sd} R_v(\mathcal{M}_{A}\oplus
\mathcal{M}_{B},r)=\min\{R_v(\mathcal{M}_{A},r),R_v(\mathcal{M}_{B},r)\}.
\end{equation}

Si en particulier, $A$ est une matrice triangulaire sup\'erieure
nilpotente \`a coefficients dans $K$, on d\'eduit du principe de
transfert de Christol ([Ch, 1.1]) que pour tout $r>0$,
\begin{equation}
\label{nl} R_v(\mathcal{M}_{(A/x)},r)=r.
\end{equation}\ \\
\textbf{Lemme 2.3.2.} \emph{Soient $n$ un entier positif, $A$ une
matrice $n\times n$ \`a coefficients dans $K(x)$. Supposons que
$\mathcal{M}_A$ est r\'egulier en $0$ avec au moins un exposant
irrationnel. Alors l'ensemble $\{p\in S\;|\;\exists \; v\in V_0,\;
v|p,\; \forall \;0<r\le 1,\;\;
R_v(\mathcal{M}_{A},r)\le r\pi_v \}$ a une densit\'e de Dirichlet positive.} \\

\emph{Preuve}. Soient $\alpha \in K$ et $r>0$. Alors la solution
$\mathcal{U}_{t_{v,r}}$ de $dU/dx=(\alpha/x) U$ au voisinage de
$t_{v,r}$, qui v\'erifie ${\mathcal U}_{t_{v,r}}(t_{v,r})=1$,  est
la s\'erie $\sum_{i\ge
0}(-1)^i(-\alpha)_i(x-t_{v,r})^i/(i!t_{v,r}^i)$.
  Donc, d'apr\`es le lemme 2.2.2, on a
\begin{equation}
\label{nr}R_v(\mathcal{M}_{(\alpha/x)},r)=
\begin{cases}
r\min(1,\pi_v\disp\liminf_{i\to+\infty} |(-\alpha)_{i}|_{v}^{-1/i})
&\mbox{si
$\alpha \in K\setminus{\bZ}_{\ge 0}$,}\\
r&\mbox{si $\alpha \in {\bZ}_{\ge 0}$.}
  \end{cases}
\end{equation}
Plus g\'en\'eralement, $R_v(\mathcal{M}_{(\alpha/x)},r)\ge r\pi_v,$
pour presque tout $ v$. Cela montre, en vertu de \eqref{tns} et
\eqref{nl}, que pour toute matrice $N$ triangulaire sup\'erieure
nilpotente \`a coefficients dans $K$, on a
\begin{equation}
\label{gen}R_v(\mathcal{M}_{\alpha/x}\otimes\mathcal{M}_{N/x},r)\ge
r\pi_v,\;\; \forall\; r>0,\;\; \text{pour presque tout}\;\; v.
\end{equation}
Par ailleurs, il est clair que $\mathcal{M}_{\alpha/x}$ n'a que des
singularit\'es r\'eguli\`eres: $0$ et l'infini. En plus, $\alpha$
est un exposant de $\mathcal{M}_{\alpha/x}$ en $0$. Donc, si $\alpha
\in K\setminus {\bQ}$ et $r=1$, le th\'eor\`eme 2 de [BS], joint \`a
\eqref{nr}, montre que l'ensemble $$\{p\in S\;|\;\exists \; v\in
V_0,\; v|p,\; \disp\liminf_{i\to\infty}
|(-\alpha)_{i}|_{v}^{-1/i}=1\}$$ a  une densit\'e de Dirichlet
positive. Ce qui entra\^{i}ne, d'apr\`es \eqref{nr}, que l'ensemble
$\{p\in S\;|\;\exists \; v\in V_0,\; v|p,\; \forall \;r>0,
R_v(\mathcal{M}_{\alpha/x},r)=r\pi_v\}$ a aussi une densit\'e de
Dirichlet positive. En combinant ceci avec \eqref{nl}, \eqref{tns},
on trouve que, pour toute matrice $N$ triangulaire sup\'erieure
nilpotente \`a coefficients dans $K$ et tout $\alpha \in K\setminus
{\bQ}$, l'ensemble $\{p\in S\;|\;\; \exists\; v\in V_0,\; v|p,\;
\forall \;r>0,\;
R_v(\mathcal{M}_{\alpha/x}\otimes\mathcal{M}_{N/x},r)=r\pi_v\}$ a
une densit\'e de Dirichlet positive. \\

Cette observation, jointe \`a \eqref{gen} et \`a \eqref{sd}, permet
d'affirmer que, pour toute matrice carr\'ee $J$ sous forme de
Jordan, \`a coefficients dans $K$, ayant au moins une valeur propre
irrationnelle, l'ensemble $\{p\in S\;|\;\;  \exists \; v\in
V_0,\;v|p,\;\;
\forall \;r>0,\;\; R_v(\mathcal{M}_{J/x},r)= r\pi_v\}$ a une densit\'e de Dirichlet positive.\\

Dans le cas g\'en\'eral, soit $A$ une matrice $n\times n$ \`a
coefficients dans $K(x)$. Supposons que $\mathcal{M}_A$ est
r\'egulier en $0$ avec au moins un exposant irrationnel. Donc,
d'apr\`es \S2.2, il existe une matrice carr\'ee $Y_A\in
{\GL}_n(K((x)))$, telle que $Y_A[A]$ est sous forme de Jordan
  et telle que $R_v(Y_A)$ est strictement positif pour tout $v\in V_0$.
  Ceci implique, d'apr\`es ce qui
pr\'ec\`ede joint \`a \eqref{chg} et le fait que $A\in
{\M}_n(\mathcal{H}_{v}(]0,1[))$ pour presque tout $v$, que
l'ensemble $\{p\in S\;|\;\; \exists \; v\in V_0,\;v|p,\; \forall
\;0<r<R_v(Y_A),\;\; R_v(\mathcal{M}_{A},r)= r\pi_v\}$  a une
densit\'e de Dirichlet positive. Par suite, le lemme r\'esulte du
fait que la fonction $\rho\mapsto
\log(R_v(\mathcal{M}_A,\exp(\rho)))$ est concave sur $]0,1]$ pour
presque tout $v$.\qed \ \\ \ \\
\textbf{Corollaire 2.3.3.} \emph{Soit $\psi\in K(x)[d/dx]$ un
op\'erateur diff\'erentiel r\'egulier en $0$. Supposons qu'il existe
une famille de  nombres r\'eels $\{r_v\}_{v\in V_0}$ de $]0,1]$
telle que $\prod_{v\in V_0}R_v(\psi,r_v)r_v^{-1}\ne 0$. Alors les
exposants  de $\psi$ en $0$ sont tous rationnels.}\\

\emph{Preuve.} Ce corollaire r\'esulte du lemme 2.3.2  et du fait
que pour tout sous-ensemble $S'\subseteq S$ non vide de densit\'e
de Dirichlet positive on a $\prod_{ p\in S'}p^{-1/(p-1)}=0$.\qed\ \\
\
\\
{\bf 2.4 Principe de transfert g\'en\'eralis\'e}\\

Le r\'esultat ci-dessous fournit, dans le cas non soluble, une
comparaison du rayon de convergence g\'en\'erique et du rayon de
convergence d'une matrice de r\'eduction  en un point r\'egulier
(cf. $[$MR, 2$]$): \\ \ \\{\bf Th\'eor\`eme 2.4.1.}\emph{ Soient
$v\in V_{0}$, $p=p(v)$ et $\phi$ un op\'erateur diff\'erentiel de
$K(x)[d/dx]$ d'ordre $n$, r\'egulier en $0$ avec exposants dans
$\mathbb{Z}_{p}$ dont les diff\'erences ne sont pas des nombres de
Liouville. Soient $Y_\phi\in {\GL}_n(K[[x]])$ une matrice de
r\'eduction de $\phi$ en $0$, $r>0$ un nombre r\'eel strictement
positif, et  $\beta_{v}$ la plus petite pente du graphe ${\mathcal
P}_{v}$ au-dessus de l'intervalle $]-\infty,\log r]$. Supposons que
$\phi$ n'ait pas de singularit\'e dans le disque \'epoint\'e
$D_{v}(0,r^{-})\setminus\{0\}$ de ${\bC}_{p}$ et que  $R
_{v}(Y_\phi)<r$. Alors}  $$(R_{v}(Y_\phi)r^{-1})^{1/n}\geq
R_{v}(\phi,r)r^{-1}\geq (R_{v}(Y_\phi)r^{-1})^{1-\beta_{v}}.$$

Notons que, sous ces conditions, le principe de transfert de
Christol  signifie que
$R_{v}(Y_\phi)=\mbox{sup}\{0<\rho<r\;|\;R_{v}(\phi,\rho)=\rho\}$.
D'autre part, puisque $Y_\phi$ est
  une matrice de r\'eduction de $\phi$, alors $x^{m}Y_\phi$ l'est aussi pour
tout entier $m\in \mathbb{Z}$. Donc on peut remplacer
  la condition $Y_\phi\in {\GL}_n(K[[x]])$  du th\'eor\`eme par $Y_\phi\in
{\GL}_n(K((x)))$. Puisque $\phi$ n'a pas de singularit\'e dans
$D_{v}(0,r^{-})\setminus\{0\}$ de ${\bC}_{p}$, le polygone
${\mathcal P}_{v}$ est alors concave sur $]-\infty,\log r]$. Donc on
a $\beta_{v}<1$, ce qui est cons\'equent avec
la borne inf\'erieure du th\'eor\`eme.\\

En combinant les th\'eor\`emes 2.3.1 et 2.4.1, on  obtient:\\
\ \\{\bf Th\'eor\`eme 2.4.2.} \emph{ Soient $\phi$ un op\'erateur de
$K(x)[d/dx]$ d'ordre $n$, r\'egulier en $0$ \`a exposants rationnels
et $(r_{v})_{v\in V_{0}}$ une famille de nombres r\'eels de $]0,1]$.
Soit $Y_\phi\in {\GL}_n(K((x)))$ une matrice de r\'eduction de
$\phi$ en $0$. Alors $$ \disp\prod_{v\in
V_{0}}\min(R_{v}(Y_\phi)r_{v}^{-1},1)\ne 0\;\;\;\;\;
\Longleftrightarrow\;\;\;\;\;\disp\prod_{v\in
V_{0}}R_{v}(\phi,r_{v})r_{v}^{-1}\ne 0.$$}\\ \emph{ Preuve.} Soit
$V_{1}$ l'ensemble de toutes les places $v\in V_{0}$ pour lesquelles
d'une part,  tous les exposants en $0$ sont dans $\mathbb{Z}_{p(v)}$
et d'autre part les singularit\'es finies non nulles ont la valeur
absolue $v$-adique 1. Soit $V_{2}$ l'ensemble de toutes les places
$v\in V_{1}$ pour lesquelles $R_{v}(Y_\phi)\geq r_{v}$,
  et soit  $\beta$ la borne inf\'erieure de l'ensemble des pentes des
polygones de convergence de $\phi$
   associ\'es aux diff\'erentes places de $V_{0}$ (dans ce cas $\beta<1$).
Donc n\'ecessairement, $V_{0}-V_{1}$ est un ensemble fini, $\phi$
n'a
    pas de singularit\'es dans le disque ouvert \'epoint\'e
$D(0,r_{v}^{-})\setminus \{0\}$ de ${\bC}_{p(v)}$ pour tout $v\in
V_{1}$ et $\phi$ v\'erifie le principe de transfert Christol pour
tout $v\in V_{2}$. Ce qui entra\^{i}ne
     $R_{v}({\phi} ,r_{v})= r_{v}$ pour tout $v\in V_{0}$. Par ailleurs,
d'apr\`es le th\'eor\`eme 2.4.1,
     on a, pour tout $v\in V_{1}-V_{2}$,
     $$(R_{v}(Y_\phi)r_{v}^{-1})^{1/n}\geq R_{v}({\phi} ,r_{v})r_{v}^{-1}
\geq (R_{v}(Y_\phi)r_{v}^{-1})^{1-\beta}.$$ Donc
$$\prod _{v\in V_{1}}R_{v}({\phi} ,r_{v})r_{v}^{-1}=
\prod _{v\in V_{1}-V_{2}}R_{v}({\phi} ,r_{v})r_{v}^{-1} \geq \prod
_{v\in V_{1}-V_{2}}(R_{v}(Y_\phi)r_{v}^{-1})^{1-\beta}=\prod _{v\in
V_{1}}\min(R_{v}(Y_\phi)r_{v}^{-1},1)^{1-\beta},$$ et
$$\prod _{v\in V_{1}}\min(R_{v}(Y_\phi)r_{v}^{-1},1)= \prod _{v\in
V_{1}-V_{2}}R_{v}(Y_\phi)r_{v}^{-1}\geq \prod _{v\in
V_{1}-V_{2}}(R_{v}({\phi} ,r_{v})r_{v}^{-1})^{n}=\prod _{v\in
V_{1}}(R_{v}({\phi} ,r_{v})r_{v}^{-1})^{n}.$$ Par cons\'equent, le
th\'eor\`eme r\'esulte du fait que $V_{0}-V_{1}$ est fini et que,
pour tout $v\in V_{0}$, on a $R_{v}(Y_\phi)>0$ et $R_{v}({\phi}
,r_{v})>0$.\qed \ \\ \ \\
{\bf 2.5 Polygone de Newton-Ramis}[Ra], [Ma, V.1]\\

Soit $\phi=\disp\sum_{i=0}^{m}
a_{i}(x)\Big(\disp\frac{d}{dx}\Big)^{i}= \disp\sum_{i=0}^{m}
\disp\sum_{j=0}^{n}a_{i,j}x^{j}\Big(\disp\frac{d}{dx}\Big)^{i}\in
K[x,\disp\frac{d}{dx}]$ un op\'erateur diff\'erentiel d'ordre $m$.
Le polygone de Newton au sens de Ramis de $\phi$ , qu'on notera
$NR(\phi)$, est l'enveloppe convexe, dans le plan $uv$, des
demi-droites horizontales $$\{u\leq i,\;v=j-i\;|\;a_{i,j}\ne 0\}.$$

A partir de cette d\'efinition, on constate que le c\^ot\'e vertical
de $NR(\phi)$ est le segment
$[(m,ord_{x}(a_{m})-m),(m,deg(a_{m})-m)]$. D'o\`u, si $NR(\phi)$ n'a
pas de c\^ot\'e vertical (autrement dit, si $a_{m}$ est un
mon\^ome), alors $\phi$ n'a pas de singularit\'e finie non nulle.
Aussi, on observe que $NR(\phi)$ est stable par homoth\'etie
(c'est-\`a-dire, par le changement de variable $x\mapsto \tau x$
o\`u $\tau\in K\setminus \{0\}$).\\

La partie de $NR(\phi)$ situ\'ee dans le demi-plan $\{v\le
ord_{x}(a_{m})-m\}$ est le polygone de Newton classique $N(\phi)$ de
$\phi$.  Le polygone $N(\phi)$ peut \^etre d\'ecrit, \`a une
translation verticale pr\`es, en consid\'erant $\phi$ comme un
op\'erateur diff\'erentiel de $K((x))[xd/dx]$ (cf. $[$Ba, 2$]$,
$[$VS, 3.3$]$). On v\'erifie que $N(\phi)$, \`a une translation
verticale pr\`es, ne d\'epend que de ${\mathcal M}_{\phi}$.
  On le note aussi $N({\mathcal M}_{\phi})$.\\

Quant \`a la partie situ\'ee dans le demi-plan $\{v\le
deg(a_{m})-m\}$, elle correspond, \`a une translation pr\`es, au
transform\'e par sym\'etrie $(u,v)\longrightarrow (u,-v)$ du
polygone de Newton de ${\mathcal M}_\phi$ en l'infini, c'est \`a
dire de $N({\mathcal
M}_{\phi}^{\infty})$.\\

La plus grande (resp. petite) pente finie de $NR(\phi)$ est
l'\textbf{invariant de Katz} de $\phi$ en $0$ (resp. en l'infini
mais avec un
signe oppos\'e) $[$De, II.1.9$]$.\\

L'op\'erateur $\phi$ est r\'egulier en $0$ si et seulement si les
pentes finies de $NR(\phi)$ sont $\le 0$. Il est r\'egulier en
l'infini si et seulement si les pentes finies de son polygone sont
$\ge 0$ (cf.
$[$Ma, III.1$]$).\\  \ \\
{\bf Lemme 2.5.1.}  \emph{Soient $a_{0}(x),\ldots,a_{m}(x)\in K[x]$.
Les op\'erateurs diff\'erentiels
  $$\phi=\sum_{i=0}^{m}
a_{i}(x)\Big(x\disp\frac{d}{dx}\Big)^{i}\;\; \mbox{et
}\;\;\;\phi_{\alpha}=\sum_{i=0}^{m}a_{i}(x)\Big(x\disp\frac{d}{dx}-\alpha\Big)^{i}
$$
ont le m\^eme polygone de Newton-Ramis}.\\ \ \\
\emph{ Preuve.} La formule
$$\phi_{\alpha}=\sum_{i=0}^{m}a_{i}(x)\sum_{j=0}^{i}\binom{i}{j}(-\alpha)^{i-j}\Big(x\disp\frac{d}{dx}\Big)^{j}=
x^{m}a_{m}(x)\Big(\disp\frac{d}{dx}\Big)^{m}+ \;\mbox{(des termes
d'ordres inf\'erieurs)}$$ montre que $NR(\phi)$ et
$NR(\phi_{\alpha})$ ont le m\^eme c\^ot\'e vertical s'il existe. On
note aussi que  ${\mathcal M}_{\phi_{\alpha}}\simeq{\mathcal
M}_{\phi}\otimes_{K(x)}{\mathcal M}_{\alpha}$. Comme ${\mathcal
M}\otimes_{K(x)}{\mathcal M}_{\alpha}$ est r\'egulier en $0$ pour
tout module
  ${\mathcal M}$ r\'egulier en $0$ (cf. $[$De, II.1.13$]$), on en d\'eduit
gr\^ace au th\'eor\`eme de
  d\'ecomposition  $[$VS, 3.54$]$ (voir aussi $[$Ma, III.1.2, III.1.5$]$,
$[$VS, 3.55$]$) que
  les polygones de Newton $N(\phi_{\alpha})$ de $\phi_{\alpha}$ et $N(\phi)$
de
  $\phi$ sont identiques \`a une translation verticale pr\`es.
   Or ces deux polygones ont le m\^eme point extr\'emal droite de
coordonn\'ees
  $(m,ord_{x}(a_{m}(x))-m)$. Donc ces polygones co\"{\i}ncident.  De m\^eme,
en utilisant le fait que
   ${\mathcal M}_{\phi_{\alpha}}^{\infty}\simeq{\mathcal
M}_{\phi}^{\infty}\otimes_{K(x)}{\mathcal
   M}_{-\alpha}$ et que les polygones $N({\mathcal
M}_{\phi_{\alpha}}^{\infty})$ et $N({\mathcal M}_{\phi}^{\infty})$
   ont le m\^eme point extr\'emal droite de coordonn\'ees
$(m,deg(a_{m}(x))-m)$, on d\'emontre que
  les polygones $N({\mathcal M}_{\phi_{\alpha}}^{\infty})$ et $N({\mathcal
M}_{\phi}^{\infty})$ co\"{\i}ncident. Par cons\'equent, $\phi$ et
$\phi_{\alpha}$ ont le m\^eme polygone de Newton-Ramis. \qed
  \ \\ \ \\
{\bf 2.6 Transformation de Fourier-Laplace g\'en\'eralis\'ee}\\

Dans ce paragraphe, on introduit le $K$-automorphisme
$\mathcal{F}_{\tau}$ de $K[x,d/dx]$, par rapport \`a l'\'el\'ement
$\tau$ de $K\setminus \{0\}$, d\'efinit par :
\begin{eqnarray*}
   \mathcal{F}_{\tau}\;\colon K[x,d/dx] &\longrightarrow & K[x,d/dx] \\
           x &\longmapsto & -\frac{1}{\tau}\frac{d}{dx}\\
           \frac{d}{dx}&\longmapsto & \tau x.
\end{eqnarray*}
Il est clair que
$\mathcal{F}_{\tau}=\mathrm{H}_\tau\circ\mathcal{F}_{1}=\mathcal{F}_{1}\circ\mathrm{H}_{1/\tau}$,
o\`u $\mathrm{H}_\tau$ est l'homoth\'etie d\'efinie par
$\mathrm{H}_\tau(x)=\tau x$. L'inverse de $ \mathcal{F}_{\tau}$, que
l'on note $\overline{\mathcal{F}_{\tau}}$, n'est autre que
$\mathrm{H}_{-1}\circ\mathcal{F}_{\tau}=\mathcal{F}_{\tau}\circ\mathrm{H}_{-1}
$. Pour $\tau=1$, on retrouve  la transform\'ee de Fourier-Laplace
classique $\mathcal{F}=\mathcal{F}_1$.

Soient $p$ un nombre premier fix\'e et  $\pi$ un nombre alg\'ebrique
tel que $\pi^{p-1}=-p$. Si $\pi\in K$, le $K$-automorphisme
$\mathcal{F}_{\pi}$
n'est autre que la transform\'ee de Fourier-Laplace normalis\'ee de Dwork.\\

Vu que le polygone de Newton-Ramis est stable par homoth\'etie, le
lemme 1.2 Chap.V de [Ma], peut s'\'etendre au cas de
$\mathcal{F}_{\tau}$:\\ \
\\{\bf Lemme 2.6.1.}\emph{ Soient $\phi\in
K[x,d/dx]$ et $\tau\in K\setminus \{0\}$. Alors
$NR(\mathcal{F}_{\tau}(\phi))$ est l'image de $NR(\phi)$ sous la
transformation
$$(u,v)\mapsto(u+v,-v).$$ En particulier, pour toute pente $t$ $($finie ou
infinie$)$ de $NR(\phi)$, le nombre
$-t/(1+t)$ est une pente de $NR(\mathcal{F}_{\tau}(\phi))$}.\\
\
\\{\bf Corollaire 2.6.2.} \emph{ Soit $\tau\in K\setminus \{0\}$. Un  op\'erateur $\phi\in K[x,d/dx]$ est
r\'egulier en $0$ et en  l'infini si et seulement si les pentes de
$NR(\mathcal{F}_{\tau}(\phi))$ appartiennent \`a $\{-1,0\}$. Dans ce
cas, il existe un entier $\ell\in\mathbb{Z}$ tel que
$x^{\ell}\mathcal{F}_{\tau}(\phi)$ peut s'\'ecrire sous la forme
$x^{\ell}{\mathcal
F}_\tau(\phi)=\Big(x\disp\frac{d}{dx}\Big)^{m}+a_{m-1}\Big(x\disp\frac{d}{dx}\Big)^{m-1}+...+a_{0}\in
K[x,x\disp\frac{d}{dx}]$.}\\

Cette derni\`ere remarque  d\'ecoule du fait que ${\mathcal
F}_\tau(\phi)$ est r\'egulier en $0$ et ne poss\`ede pas de
singularit\'e finie non nulle.
\section{ $E$-fonctions et $E$-op\'erateurs}
{\bf 3.1 $E$-fonctions} \\

Rappelons qu'une s\'erie de $\sum_{n\ge 0}a_nx^n\in
\overline{\bQ}[[x]]$ est dite holonome si elle est solution d'un
op\'erateur diff\'erentiel de $\overline{\bQ}[x,d/dx]$. Les
coefficients d'une telle s\'erie v\'erifient une relation de
r\'ecurrence de la forme $\sum_{0\le i\le \ell}P_i(n+i)a_{n+i}=0$,
pour tout entier $n\ge 0$ o\`u les $P_0,\ldots,P_\ell$ d\'esignent
des polynomes \`a coefficients alg\'ebriques. Cela veut dire que
tous les coefficients $a_n$
appartiennent \`a une extension finie de ${\bQ}$.\\

Dans ce qui suit, on s'interessera juste aux s\'eries Gevrey
arithm\'etiques holonomes. Pour cela, on supposera que  le corps de
nombres $K$ contienne tous les coefficients des s\'eries Gevrey
arithm\'etiques holonomes qu'on rencontrera par la suite.
 L'ensemble des s\'eries Gevrey arithm\'etiques d'ordre
 $s\in {\bQ}$ \`a coefficients dans $K$ sera not\'e $K\{x\}_s$.\\ \ \\
{\bf Lemme 3.1.1.}\emph{ {\rm (1)} Si $g$ est une $G$-fonction,
alors
$\disp\prod_{v\in V_{0}} \min(r_{v}(g),1)\ne 0$ .\\
{\rm (2)} Soit $s\in {\bQ}$. Les sommes finies $\disp\sum_{\alpha,k}
\lambda_{\alpha,k}y_{\alpha,k}x^{\alpha}\log^{k}x$, o\`u les
$\alpha\in {\mathbb{Q}}, k\in {\mathbb{N}}$, $\lambda_{\alpha,k}\in
K$ et les $y_{\alpha,k}\in K\{x\}_s$, forment une $K[x]$-alg\`ebre
diff\'erentielle que l'on note $NGA\{x\}_{s}$.}\\\ \\ (cf.
$[$A1, p126$])$, $[$A2, 1.4, 2.1$]$).\\

En combinant ce lemme avec le lemme 2.2.2, on trouve : \\
\ \\{\bf Corollaire 3.1.2.} \emph{Si $F$ est une $E$-fonction,
alors} $\disp\prod_{v\in V_{0}}
\min(r_{v}(F)\pi_{v}^{-1},1)\ne 0$.\\
\ \\
{\bf 3.2 $G$-op\'erateurs}\\

Les $G$-op\'erateurs poss\`edent aussi une caract\'erisation
$p$-adique dite ``local Bombieri'', \'equivalente \`a la condition
``Bombieri", que nous donnons ici au voisinage de l'infini:
\\ \
\\{\bf Proposition 3.2.1.} \emph{Soit $\phi$ un op\'erateur diff\'erentiel de
$K(x)[d/dx]$ d'ordre $m$,  ayant une singularit\'e r\'eguliere en
l'infini avec exposants rationnels. Alors les conditions suivantes
sont \'equivalentes}
\begin{itemize}
\item[(1)] $\prod_{v\in V_{0}}R_{v}(\phi,1)\ne 0$. \item[(2)]
$\phi$  \emph{poss\`ede une matrice de r\'eduction en
l'infini} $Y_{\phi}\in {\GL}_m(K((x)))$ \emph{telle que}\\
$\prod_{v\in V_{0}}\min(R_{v}(Y_{\phi}),1)\ne 0$.
\end{itemize}
\emph{Preuve.}  Cette proposition d\'ecoule  du th\'eor\`eme 2.4.2,
et
de la remarque suivante: \qed \ \\ \ \\
{\bf Remarque 3.2.2.} Si $\phi^{\infty}$ d\'esigne l'op\'erateur
obtenu \`a partir d'un op\'erateur $\phi$ de $K(x)[d/dx]$ par le
changement de variable $x\longrightarrow 1/x$, alors on v\'erifie
ais\'ement que  $R_{v}(\phi^{\infty},r^{-1})=R_{v}(\phi,r)r^{-2}$,
pour tout $r>0$ et tout $v\in V_{0}$. \\ \ \\
{\bf 3.3 $E$-op\'erateurs}\\

Il r\'esulte de  la d\'efinition des $E$-op\'erateurs et des
propri\'et\'es des $G$-op\'erateurs qu'un produit $\psi.\psi'$ de
$E$-op\'erateurs est un $E$-op\'erateur, que tout diviseur \`a
droite d'un $E$-op\'erateur est un $E$-op\'erateur et que toute
paire de $E$-op\'erateurs poss\`ede un multiple \`a gauche commun
qui est un $E$-op\'erateur. De m\^eme que pour les $G$-op\'erateurs,
si $\psi$ est un $E$-op\'erateur, alors $\psi^{*}$ et
$\overline{\psi}$ sont aussi des $E$-op\'erateurs (cf. $[$A2,
4.1$]$).\\
\
\\{\bf Lemme 3.3.1.} \emph{Soit $\psi\in K[x,d/dx]$ un
$E$-op\'erateur. Alors $x^{\ell}\psi$ est un $E$-op\'erateur pour
tout  entier $\ell$  tel que $x^{\ell}\psi\in K[x,d/dx]$.}\\\
\\Cela d\'ecoule du fait que, Si $\ell< 0$, L'op\'erateur $x^{\ell}\psi$ est
un diviseur \`a gauche de $\psi$, tandis que, si $\ell> 0$, cet
op\'erateur est le produit des deux $E$-op\'erateurs $x^{\ell}$
et $\psi$. \\

A partir du th\'eor\`eme  ci-dessous d\^u \`a Y. Andr\'e (cf. $[$A2,
4.3$]$), on va pouvoir d\'emontrer que les $E$-op\'erateurs
poss\`edent des propri\'et\'es analogues aux  propri\'et\'es ``local
Galochkin",
  ``Local Bombieri" et ``Bombieri"  des $G$-op\'erateurs .
\\ \ \\
{\bf Th\'eor\`eme 3.3.2}.\emph{ Soit $\psi$ un $E$-op\'erateur
d'ordre m, alors :
\begin{itemize}
\item[(1)] toutes les pentes de $NR(\psi)$ sont dans $\{-1,0\}$,
en particulier $\psi$ n'a pas de singularit\'e finie non nulle;
\item[(2)] $\psi$ est r\'egulier en $0$ avec exposants rationnels;
\item[(3)] $\psi$ admet une base de solutions en $0$ de la forme
$(F_{1},...,F_{m})x^{\Gamma_0}$ o\`u les $F_{i}$ sont des
$E$-fonctions, et $\Gamma_0$ d\'esigne une matrice carr\'ee d'ordre
$m$ triangulaire sup\'erieure \`a coefficients dans $\mathbb{Q}$;
\item[(4)] $\psi$ admet une base de solutions en l'infini de la
forme
$$\Big(f_1\Big(\frac{1}{x}\Big),\ldots,f_m\Big(\frac{1}{x}\Big)\Big)\Big(\frac{1}{x}\Big)^{\Gamma_\infty}.e^{-\Delta\; x}$$
o\`u les $f_i$ sont des \'el\'ements de $K\{x\}_1$, o\`u $\Delta$
est une matrice diagonale ayant pour coefficients diagonoaux les
singularit\'es de $\overline{\mathcal{F}}(\psi)$, et o\`u
$\Gamma_\infty$ d\'esigne une matrice carr\'ee triangulaire
sup\'erieure \`a coefficients dans $\bQ$, qui commute \`a $\Delta$.
\end{itemize}
} \

Ce th\'eor\`eme montre que tout $E$-op\'erateur poss\`ede une
base de solutions en $0$ \`a coefficients dans $NGA\{x\}_{-1}$.\\
\
\\{\bf Corollaire 3.3.3.} \emph{Soient $m\in {\bZ}_{\ge 0}$, $\alpha\in
\mathbb{Q}$ et $\psi=\disp\sum_ {i=0}^{m}
a_{i}(x)\Big(x\disp\frac{d}{dx}\Big)^{i}\in K[x,x\disp\frac{d}{dx}]$
avec $a_m(x)=1$. Alors, $\psi$ est un $E$-op\'erateur si et
seulement si $\psi_{\alpha}:=\disp\sum_{i=0}^{m}
a_{i}(x)\Big(x\disp\frac{d}{dx}-\alpha\Big)^{i}$ est un
$E$-op\'erateur.}\\ \ \\
\emph{ Preuve.} Supposons que $\psi$ est un $E$-op\'erateur.
D'apr\`es le th\'eor\`eme ci-dessus,
  $\psi$  poss\`ede une base de
solutions $(\zeta_{1},\ldots,\zeta_{m})$ \`a coefficients dans
$NGA\{x\}_{-1}$. Puisque $\psi
x^{-\alpha}=x^{-\alpha}\psi_{\alpha}$, on en d\'eduit que
$(\zeta_{1},\ldots,\zeta_{m})x^{\alpha}$ est une base de solutions
de $\psi_{\alpha}$ \`a coefficients dans $NGA\{x\}_{-1}$. Donc en
vertu de $[$A2, 6.1$]$, il existe des $E$-op\'erateurs $\psi_{i}$
tels que $\psi_{i}(\zeta_{i}x^{\alpha})=0,$ pour $i=1,\ldots,m$. De
plus, ces op\'erateurs
   admettent un multiple commun
  \`a gauche $\psi_{0}$ d'ordre $n\ge m$, qui est aussi, d'apr\`es \S3.3, un
$E$-op\'erateur.
  Par ailleurs, et gr\^ace au lemme 3.3.1, il existe un entier
$\ell\in\mathbb{Z}$ tel que $x^{\ell}\psi_{0}$
   est un $E$-op\'erateur qui appartient \`a $K[x,x\disp\frac{d}{dx}]$.
D'autre part, puisque $a_{m}(x)=1$, il existe $P,Q \in
K[x,x\disp\frac{d}{dx}] $ tels que
$x^{\ell}\psi_{0}=P\psi_{\alpha}+Q$ et
$ord(Q)<ord(\psi_{\alpha})=m$. Or, pour $i=1,\ldots,m$, on a
$\psi_{0}(\zeta_{i}x^{\alpha})=0$. On en d\'eduit que $Q=0$,
c'est-\`a-dire que $x^{\ell}\psi_{0}=P\psi_{\alpha}$. Donc
$\psi_{\alpha}$ est un diviseur \`a droite d'un $E$-op\'erateur et,
par cons\'equent, il est aussi un $E$-op\'erateur d'apr\`es \S3.3.
La r\'eciproque
d\'ecoule du fait que $(\psi_{\alpha})_{-\alpha}=\psi$.\qed\ \\
\section{Conditions n\'ecessaires}
Le th\'eor\`eme 3.3.2 montre que la condition (1) du th\'eor\`eme
1.1  est  n\'ecessaire. Dans cette section
on d\'emontrera que la deuxi\`eme condition l'est aussi.\\
\
\\{\bf Lemme 4.1.}\emph{ Soient ${\mathcal M}$ un $K(x)$-module
diff\'erentiel r\'egulier en $0$ et $A\in {\M}_m(K(x))$ une matrice
associ\'ee. Soient $Y$ et $Z$ deux matrices de r\'eduction de $A$ en
$0$ dans ${\GL}_m(K((x)))$. Alors il existe une extension finie $K'$
de $K$ et une
matrice $L\in {\GL}_m(K'[x,1/x])$ telle que $Y=LZ$. }\\
\ \\\emph{ Preuve.} Par d\'efinition, on a: $$
Y[A]=\disp\frac{1}{x}\Gamma_{1},
\;\;\;\;Z[A]=\disp\frac{1}{x}\Gamma_{2}\;\;\; \mbox{avec}\;\;\;
\Gamma_{1},\Gamma_{2}\in{\M}_m(K).$$ Posons $L=YZ^{-1}$. On en
d\'eduit,
$$L[\disp\frac{1}{x}\Gamma_{2}]=\disp\frac{1}{x}\Gamma_{1}\;\;\;\;\mbox{c'est-\`a-dire}\;\;\;\;
x\disp\frac{d}{dx}L=\Gamma_{1}L-L\Gamma_{2}.$$ Donc, si on \'ecrit
$L= \disp\sum_{i\in\mathbb{Z}}L_{i}x^{i} $, on trouve pour tout
$i\in \mathbb{Z}\setminus\{0\}$
$$iL_{i}=\Gamma_{1}L_{i}-L_{i}\Gamma_{2}.$$ Or les valeurs propres
de l'application
\begin{eqnarray*}
   T_{i}\;\colon {\M}_m(K) &\longrightarrow &{\M}_m(K) \\
           X &\longmapsto &\Gamma_{1}X-X\Gamma_{2}-iX
\end{eqnarray*}
sont $\{\lambda-\gamma-i\}$ (o\`u les $\lambda$ et $\gamma$ sont
respectivement des valeurs propres de $\Gamma_{1}$ et de
$\Gamma_{2}$). Donc, sauf peut-\^etre pour un nombre fini d'entiers
$i$, cette application est inversible. Par cons\'equent, $L_{i}$ est
nul \`a l'exception d'un
nombre fini d'entiers $i$. La conclusion suit.\qed \ \\

On pourra toujours supposer que $K$ est assez grand de sorte qu'il
contienne les coefficients d'une \'eventuelle matrice $L$ comme
ci-dessus.
\\ \ \\
{\bf Lemme 4.2.}\emph{ Soient $A$ et $B$ deux matrices de
${\M}_m(K(x))$ associ\'ees \`a un $K(x)$-module diff\'erentiel
$\mathcal{M}$, et soient $Y_A$ et $Y_B$ respectivement deux matrices
de r\'eductions de $A$ et $B$ en $0$. Alors,
\begin{eqnarray}
\min(R_v(Y_A),1)&=&\min(R_v(Y_B),1)
\end{eqnarray}
pour presque tout  $v$ de $V_0$.}
\\ \ \\
\emph{Preuve.} Par hypoth\`ese, il existe une matrice $H$ de
${\GL}_m(K(x))$ telle que $H[A]=B$. Donc $Y_BH$ est aussi une
matrice de r\'eduction de $A$ en $0$. D'apr\`es le lemme 4.1, il
existe une matrice $L$ de ${GL}_m(K[x,1/x])$ telle que
\begin{eqnarray}
Y_A=LY_BH.
\end{eqnarray}
Par ailleurs, les coefficients de $H$ et $H^{-1}$, consid\'er\'es
comme \'etant des \'el\'ements de $K((x))$,
  ont des rayons de convergence $\ge 1$ pour presque
tout $v$ de $V_0$. Donc, en vertu de (4.2) on a,
  $$r_v(Y_A)\ge\min(r_v(Y_B),r_v(H),r_v(L))\ge\min(r_v(Y_B),1)$$
  et  aussi $$r_v(Y_A^{-1})\ge\min(r_v(Y_B^{-1}),1),$$
  pour presque tout $v$ de $V_0$. Ce qui
donne   $$\min(R_v(Y_A),1)\ge
   \min(R_v(Y_B),1)$$
pour presque tout $v$ de $V_0$. L'in\'egalit\'e dans l'autre sens
s'obtient du fait que $ Y_B=L^{-1}Y_AH^{-1}$ et la conclusion
suit.\qed\
\\\ \\
\textbf{Corollaire 4.3.} \emph{Sous les hypoth\`eses du lemme 4.2,
on a les \'equivalences suivantes:}
\begin{eqnarray}
\prod_{v\in V_0}\min(R_v(Y_A),1)\ne 0&\Longleftrightarrow&
\prod_{v\in
V_0}\min(R_v(Y_B),1)\ne 0\\
\prod_{v\in V_0}\min(R_v(Y_A)\pi_v^{-1},1)\ne
0&\Longleftrightarrow&\prod_{v\in V_0}\min(R_v(Y_B)\pi_v^{-1},1)\ne
0
\end{eqnarray}
\emph{Preuve} Ces \'equivalences d\'ecoulent de (4.1) et de la
proposition 2.2.1.\qed \ \\ \ \\
{\bf Lemme 4.4.} \emph{Soit $\psi$ un op\'erateur diff\'erentiel de
$K[x,d/dx]$ d'ordre $m$, et soient $(y_1,\ldots,y_m)x^{\Gamma_0}$ et
$(z_1,\ldots,z_m)x^{\Lambda_0}$ respectivement des bases de
solutions de $\psi$ et $\psi^*$ en $0$, o\`u les coefficients
$y_1,\ldots,y_m,z_1,\ldots,z_m$ appartiennent \`a $K[[x]]$ et o\`u
$\Gamma_0$ et $\Lambda_0$ sont des  matrices $m\times m$ \`a
coefficients dans $K$. Alors $\psi$ poss\`ede une matrice de
r\'eduction $Y_\psi\in{\GL}_m(K[[x]])$ en $0$ telle que les
coefficients de $Y_\psi$ $($resp. de $Y_\psi^{-1})$ soient des
combinaisons $K[x,1/x]$-lin\'eaires de $z_1,\ldots,z_m$ $($resp.
$K[x]$-lin\'eaires de $y_1,\ldots,y_m)$ et de
leurs d\'eriv\'ees.}\\
\emph{Preuve.} \'Ecrivons
$\psi=\disp\sum_{i=0}^{m}a_{i}(x)\Big(\disp\frac{d}{dx}\Big)^{i}\in
K[x,\disp\frac{d}{dx}]$. Comme  $(y_1,\ldots,y_m)x^{\Gamma_0}$ est
une base de solutions  de $\psi$ en $0$, la matrice wronskienne $W$
de $(y_1,\ldots,y_m)x^{\Gamma_0}$ est solution du syst\`eme
$dX/dx=A_{\psi}X$. De plus, il existe un entier $m$ pour lequel $W$
peut s'\'ecrire sous la forme $Y_{0}x^{\Gamma}$ o\`u
$\Gamma=\Gamma_0+m{\bI}_m$ et  o\`u $Y_{0}$ est une matrice de
${\GL}_m(K[[x]])$ dont les coefficients sont des combinaisons
$K[x]$-lin\'eaires de  $y_1,\ldots,y_m$ et de leurs d\'eriv\'ees.
Alors on a, $Y_{0}^{-1}[A_\psi]=\Gamma/x$, ce qui veut dire que
$Y_{0}^{-1}$ est une matrice de r\'eduction de $\psi$ en $0$, ou
encore que
\begin{eqnarray}
^{t}Y_{0}[-^{t}A_{\psi}]&=&-^{t}\Gamma\disp\frac{1}{x}.
\end{eqnarray}
Par ailleurs, le $m$-uplet  $a_{m}(x) (z_1,\ldots,z_m)x^{\Lambda_0}$
est une base de solutions de
$\psi^{*}a_{m}^{-1}=(a_{m}^{-1}\psi)^{*}$ en $0$. Donc, la matrice
$U$ dont les lignes $u_1,\ldots,u_m$ sont d\'efinies par les
relations
\begin{eqnarray*}
u_m&=&a_{m}(x)
(z_1,\ldots,z_m)x^{\Lambda_0}\\
u_{m-i}&=&\disp\frac{a_{m-i}(x)}{a_{m}(x)}u_{m}-\disp\frac{d}{dx}u_{m-i+1}\;\;\;(1\le
i\le m-1),
\end{eqnarray*}
est une solution du syst\`eme $dX/dx=-^{t}A_{\psi}X.$ De plus, $U$
peut s'\'ecrire sous la forme $Zx^{\Lambda_0}$, o\`u $Z$ est une
matrice  inversible $m\times m$  dont les coefficients sont des
combinaisons $K[x,1/x]$-lin\'eaires de $z_1,\ldots,z_m$ et de leurs
d\'eriv\'ees. Donc, on a
\begin{eqnarray}
Z^{-1}[-^{t}A_{\psi}]&=&\disp\frac{1}{x}\Lambda_0.
\end{eqnarray}
Par cons\'equent, d'apr\`es les formules (4.5), (4.6) et le lemme
4.1, il existe une matrice $L\in {\GL}_m(K[x,1/x])$ telle que
$^{t}Y_{0}=LZ^{-1}$. Par suite, les \'el\'ements de $Y_{0}^{-1}$
sont des  combinaisons $K[x,1/x]$-lin\'eaires de $z_1,\ldots,z_m$ et
de leurs d\'eriv\'ees. La matrice $Y_\psi:=Y_{0}^{-1}$ poss\`ede
donc toutes les propri\'et\'es
annonc\'ees.\qed  \ \\ \ \\
{\bf Th\'eor\`eme 4.5.}\emph{ Soit $\psi$ un $E$-op\'erateur de
$K[x,d/dx]$ d'ordre $m$, alors :
\begin{itemize}
\item[(1)] $\psi$ poss\`ede une matrice de r\'eduction $Y_\psi\in
{\GL}_m(K((x)))$ en $0$ telle que les coefficients de $Y_\psi$ et
$Y_\psi^{-1}$ sont des $E$-fonctions; \item[(2)] $\prod_{v\in V_{0}}
\min( R_{v}(Y_\psi)\pi_{v}^{-1},1)\ne 0;$ \item[(3)] $\prod_{v\in
V_{0}} R_{v}(\psi,\pi_{v})\pi_{v}^{-1}\ne 0. $
\end{itemize}
}\ \\ \emph{ Preuve.} (1) r\'esulte du th\'eor\`eme 3.3.2 et des
lemmes 3.1.1 et 4.4.\\(2) \'Ecrivons $Y_\psi=(y_{ij})_{ij}$ et
$Y_\psi^{-1}=(\widetilde{y}_{kl})_{kl}$.  On a
\begin{eqnarray*}
\min(R_{v}(Y_\psi)\pi_{v}^{-1},1)&=&\min(\min_{i,j}(r_{v}(y_{ij})\pi_{v}^{-1}),\min_{k,l}(r_{v}(\widetilde{y}_{kl})\pi_{v}^{-1}),1)\\
                 &\ge
&\prod_{ij}\min(r_{v}(y_{ij})\pi_{v}^{-1},1)\prod_{kl}\min(r_{v}(\widetilde{y}_{kl})\pi_{v}^{-1},1)
\end{eqnarray*}
et, gr\^{a}ce \`a (1) et le corollaire 3.1.2, on obtient
$$\prod_{v\in
V_{0}}\min(R_{v}(Y_\psi)\pi_{v}^{-1},1)\ge\prod_{ij}\prod_{v\in
V_{0}}\min(r_{v}(y_{ij})\pi_{v}^{-1},1)\prod_{kl}\prod_{v\in
V_{0}}\min(r_{v}(\widetilde{y}_{kl})\pi_{v}^{-1},1)\ne 0.
$$
L'assertion (3) se d\'eduit de (2) et du th\'eor\`eme 2.4.2.\qed\ \\

Un op\'erateur $\psi$ de $K[x,d/dx]$ qui n'a pas de singularit\'e
finie non nulle et qui remplit la condition (2) du th\'eor\`eme
ci-dessus n'est pas n\'ecessairement un $E$-op\'erateur comme le
montre l'exemple
suivant :\\ \ \\
{\bf Exemple 4.6. } Soit $\psi =xd/dx-a-bx^{2}$ avec $a,b\in
\mathbb{Q}-\{0\}$. D'une part, $\psi$ poss\`ede seulement deux
singularit\'es : le point $0$, qui est r\'egulier d'exposant $a$, et
l'infini qui est irr\'egulier. D'autre part, la fonction
$Y=\exp(bx^{2}/2)$ est une matrice de  r\'eduction de $\psi$ en $0$,
puisque $Y\Big[(a+bx^{2})/x\Big]=a/x$. De plus, on trouve
$R_{v}(Y)=(|2b^{-1}|_{v}\pi_{v})^{1/2}$ pour toute place $v$ de
$V_{0}$.  Donc, d'apr\`es le th\'eor\`eme 2.4.2, on a $ \prod_{v\in
V_{0}}R_{v}(\psi,\pi_{v})\pi_{v}^{-1}\ne 0.$ Cependant, l'invariant
de Katz de $\psi$ en l'infini est \'egal \`a $2$. Donc, $-2$ est une
pente
   de $NR(\psi)$. Par suite, $\psi$
    n'est pas un $E$-op\'erateur.
\section{Transformation de Laplace formelle g\'en\'eralis\'ee}
\label{sec:condsuff} Fixons d'abord un plongement de $K$ dans $\bC$
et un nombre $\lambda\in K\setminus\bZ$. Pour tout $i\in\bZ$ tel que
la partie r\'eelle de $\lambda+i$ soit $>-1$, la transform\'ee de
Laplace (au sens complexe) du mon\^ome $x^{\lambda+i}$ est donn\'ee
par
\[
\mathcal{L}(x^{\lambda +i})
 = \Gamma(\lambda+i+1)x^{-\lambda -i-1}
 = \Gamma(\lambda)(\lambda)_ix^{-\lambda -i-1},
\]
o\`u $(\lambda)_i = \lambda(\lambda+1)\cdots (\lambda+i-1)$ est un
\'el\'ement de $K\setminus\{0\}$.  Le produit $\Gamma(\lambda)^{-1}
\mathcal{L}$ s'\'etend donc naturellement en un op\'erateur de
$x^{\lambda}K[[x]]$ dans $x^{-\lambda-1}K[[1/x]]$.  Plus
g\'en\'eralement, \'etant donn\'e $\lambda\in K\setminus\bZ$ et
$\tau\in \overline{\bQ}\setminus \{0\}$, on d\'efinit
\begin{eqnarray*}
   \mathcal{L}_{\lambda}^{\tau}\;\colon x^{\lambda}K[[x]]
   &\longrightarrow & \frac{1}{x^{\lambda+1}}K(\tau)[[1/x]] \\
           \sum_{i\ge 0}a_ix^{\lambda +i} &\longmapsto &
           \sum_{i\ge 0}a_i
\frac{(\lambda)_i}{\tau^{i}}\frac{1}{x^{\lambda+i+1}}.
\end{eqnarray*}
Cet op\'erateur satisfait les relations, de commutation avec la
d\'erivation $d/dx$, suivantes: $$\forall\;\; f \in
x^{\lambda+1}K[[x]],\;\;\;\;\;{\mathcal
L}^{\tau}_\lambda\Big(\frac{df}{dx}\Big) = {\tau}x {\mathcal
L}^{\tau}_\lambda(f)\;\;\;\;\;\text{et}\;\;\;\;\;{\mathcal
L}^{\tau}_\lambda(xf) = -\frac{1}{{\tau}}\frac{d}{dx} {\mathcal
L}^{\tau}_\lambda(f).$$ Dans le cas ou $\tau=1$, on retrouve les
relations de commutation, bien connues, de la transform\'ee de
Laplace avec $d/dx$.

Pour un nombre premier $p$ fix\'e et un nombre alg\'ebrique $\pi\in
\overline{\bQ}$ tel que $\pi^{p-1}=-p$, on constate
 que l'op\'erateur ${\mathcal L}^{-\pi}_\lambda$ est diff\'erent de la
transformation de Laplace ``modifi\'ee" de Dwork dans le cas d'une
seule variable [Dw, 11.3].

Plus g\'en\'eralement, fixons $\tau\in K\setminus\{0\}$ et une
matrice carr\'ee $\Lambda$ de format $n\times n$ dont les valeurs
propres appartiennent \`a $K\setminus {\bZ}$.  Dans cette section,
on g\'en\'eralise cette construction \`a des produits
$A(x)x^\Lambda$ o\`u $A(x)$ est une matrice $m\times n$ \`a
coefficients dans $K[[x,1/x]]$.  On obtient ainsi une transformation
de Laplace formelle
$$
 \mathcal{L}_{\Lambda}^{\tau}\; :\; x^{\Lambda}{\M}_{m\times
n}(K[[x,1/x]])\;\longrightarrow\;
\frac{1}{x^{\Lambda+1}}{\M}_{m\times n}(K[[x,1/x]])
$$
qui poss\`ede des propri\'et\'es de commutation avec la d\'erivation
$d/dx$ (5.4) \'etendant celles du cas ($m=n=1$) ci-dessus, et qui
pr\'eserve par suite la dualit\'e entre la transform\'ee de Laplace
des fonctions et celle de Fourier-Laplace des op\'erateurs
(Proposition 5.1). On compare aussi la couronne
 de convergence de la matrice $A(x)$ et celle de la partie
analytique de $\mathcal{L}_{\Lambda}^{\tau}(A(x)x^\Lambda)$
(Proposition 5.2.2). Dans le cas o\`u $m=1$ et $A(x)x^\Lambda$ est
une base de solutions logarithmiques d'un op\'erateur diff\'erentiel
donn\'e $\phi \in K[x,d/dx]$  en $0$ (resp. en l'infini), les
coefficients de $\mathcal{L}_{\Lambda}^{\tau}(A(x)x^\Lambda)$ sont
des solutions logarithmiques de $\mathcal{F}_{\tau}(\phi)$ en
l'infini (resp. en $0$). Ceci permet d'avoir plus d'informations sur
$\mathcal{F}_{\tau}(\phi)$ (exposants, rayon de convergence des
coefficients des solutions logarithmiques) \`a partir de celles de
$\phi$ et vice-versa. Dans le cas particulier o\`u $\tau=1$, on
observe que cette transformation joue presque le m\^eme r\^ole que
celui de la transform\'ee de Laplace \'etendue \`a l'aide du calcul
op\'erationnel [A2, 5.3]. L'avantage de cette nouvelle
transformation est qu'elle ne fait pas intervenir des coefficients
transcendants.\\

Si le nombre $\pi$, d\'efini ci-dessus, appartient \`a $K$, il nous
semble que cette nouvelle transformation sera utile pour l'\'etude
des propri\'et\'es $p$-adiques de la transform\'ee de
Fourier-Laplace normalis\'ee de Dwork $\mathcal{F}_\pi$.\\

Enfin, cette transformation convient parfaitement \`a notre
situation pour compl\'eter la preuve du th\'eor\`eme principal de
cet article (voir section 6).

\subsection{5.1 Transformation de Laplace formelle}\ \\

Soient $n$ un entier positif, $\tau$ un \'el\'ement de
$K\setminus\{0\}$ et $\Lambda$ une matrice $n\times n$ \`a
coefficients dans $\bQ$ dont toutes les valeurs propres
appartiennent \`a $\bQ\setminus \bZ$ (la transformation qu'on va
introduire s'\'etend d'une mani\`ere naturelle au cas o\`u
$\Lambda\in {\GL}_n(K)$ avec des valeurs propres dans $K\setminus
{\bZ}$  tout en pr\'eservant les propri\'et\'es de la proposition
5.1 ci-dessous). Soit $d\ge 1$ un d\'enominateur commun des valeurs
propres de $\Lambda$. Alors $x^\Lambda$ est une matrice de
${\GL}_n(\bQ(x^{1/d},\log x))$ qui satisfait
\begin{equation*}
\frac{d}{dx}(x^\Lambda) = \Lambda x^{-1}x^{\Lambda} = \Lambda
x^{\Lambda-{\bI}_n}
\end{equation*}
(voir \S8, chap.\ III de \cite{DGS}). Pour tout entier $m\ge 1$ et
toute matrice $A(x) = \sum_{i=-\infty}^\infty A_i x^i$ de format
$m\times n$ \`a coefficients dans $K[[x,1/x]]$, on d\'efinit la
\emph{transform\'ee de Laplace de $A(x)x^\Lambda$ par rapport \`a}
$\Lambda$ et \`a ${\tau}$ par
\begin{equation}
\label{trf} {\mathcal L}^{\tau}_\Lambda(A(x)x^\Lambda) = {\mathcal
L}^{\tau}_\Lambda\Big(\sum_{i=-\infty}^\infty A_i
x^{\Lambda+i{\bI}_n}\Big) := \sum_{i=-\infty}^\infty A_i
C_{\Lambda,\tau}(i) x^{-\Lambda-(i+1){\bI}_n}
\end{equation}
o\`u $C_{\Lambda,\tau}\colon\bZ\to\GL_n(\bQ(\tau))$ est la fonction
d\'etermin\'ee par les conditions
\begin{equation}
\label{def:Cgamma}C_{\Lambda,\tau}(0)={\bI}_n \et
C_{\Lambda,\tau}(i)={\tau}^{-1}(\Lambda+i{\bI}_n)C_{\Lambda,\tau}(i-1),
\quad (i\in\bZ).
\end{equation}
Explicitement, on a
\begin{equation*}
C_{\Lambda,\tau}(i) =
  \begin{cases}
  \tau^{-i} (\Lambda+i{\bI}_n)(\Lambda+(i-1){\bI}_n)\cdots (\Lambda+{\bI}_n) &\mbox{si
$i\ge 1$,}\\
   {\bI}_n &\mbox{si $i=0$,}\\
   \tau^{-i} \Lambda^{-1}(\Lambda-{\bI}_n)^{-1}\cdots(\Lambda+(i+1){\bI}_n)^{-1}
&\mbox{si
   $i\le -1$.}
  \end{cases}
\end{equation*}
Gr\^ace \`a ces formules, ou par r\'ecurrence sur $i$ sur la base de
\eqref{def:Cgamma}, on observe que $C_{\Lambda,\tau}$  prend ses
valeurs dans un sous-anneau commutatif de $\GL_n(\bQ(\tau))$ et que
\begin{equation}
\label{CC} C_{\Lambda,\tau}(i) C_{-\Lambda,\tau}(-i-1) =
(-1)^{i+1}\tau\Lambda^{-1}
\end{equation}
pour tout $i\in\bZ$. La transformation ${\mathcal L}^{\tau}_\Lambda$
poss\`ede les propri\'et\'es formelles suivantes, analogues \`a
celle de la transformation de Laplace usuelle.\\ \ \\
{\bf Proposition 5.1.1.} \emph{Soient $m$ un entier positif et
$A(x)$ une matrice $m\times n$ \`a coefficients dans $K[[x,1/x]]$.
On pose $f=A(x) x^{\Lambda}$. Alors, on a
\begin{equation}
\label{Ld=xL} {\mathcal L}^{\tau}_{-\Lambda}\big({\mathcal
L}^{\tau}_\Lambda(f)\big) = - \tau A(-x)\Lambda^{-1} x^\Lambda,
\quad {\mathcal L}^{\tau}_\Lambda\Big(\frac{df}{dx}\Big) = {\tau}x
{\mathcal L}^{\tau}_\Lambda(f)\;\; \text{et}\;\; {\mathcal
L}^{\tau}_\Lambda(xf) = -\frac{1}{{\tau}}\frac{d}{dx} {\mathcal
L}^{\tau}_\Lambda(f).
\end{equation}
De plus, si $m=1$ et si les coefficients de $f$ sont des solutions
d'un op\'erateur diff\'erentiel $\phi\in K[x,d/dx]$, alors ceux de
${\mathcal L}^{\tau}_\Lambda(f)$ sont des solutions de ${\mathcal
F}_{\tau}(\phi)$.}
\begin{proof}[Preuve] Par lin\'earit\'e, il suffit de d\'emontrer
\eqref{Ld=xL} dans le cas o\`u $f=x^ix^\Lambda=x^{\Lambda+i{\bI}_n}$
avec $i\in\bZ$.  Dans ce cas, gr\^ace \`a \eqref{def:Cgamma} et
\eqref{CC}, on trouve bien
\begin{align*}
{\mathcal L}^{\tau}_{-\Lambda}\big({\mathcal L}_\Lambda(f)\big)
  &= {\mathcal L}^{\tau}_{-\Lambda}\big(C_{\Lambda,\tau}(i)x^{-\Lambda-(i-1){\bI}_n}\big)
  = C_{\Lambda,\tau}(i)C_{-\Lambda,\tau}(-i-1)x^{\Lambda+i{\bI}_n}
  =- \tau A(-x)\Lambda^{-1} x^\Lambda,\\
{\mathcal L}^{\tau}_\Lambda\Big(\frac{df}{dx}\Big)
  &= (\Lambda+i{\bI}_n) C_{\Lambda,\tau}(i-1) x^{-\Lambda-i{\bI}_n}
  = {\tau}C_{\Lambda,\tau}(i) x^{-\Lambda-i{\bI}_n}
  = {\tau}x {\mathcal L}^{\tau}_\Lambda(f),\\
{\mathcal L}^{\tau}_\Lambda(xf) &= C_{\Lambda,\tau}(i+1)
x^{-\Lambda-(i+2){\bI}_n}
  = \frac{1}{{\tau}}C_{\Lambda,\tau}(i) (\Lambda+(i+1){\bI}_n) x^{-\Lambda-(i+2){\bI}_n}
  = -\frac{1}{{\tau}}\frac{df}{dx} {\mathcal L}^{\tau}_\Lambda(f).
\end{align*}
Enfin, les formules \eqref{Ld=xL} impliquent ${\mathcal
L}^{\tau}_\Lambda(\phi(f)) = {\mathcal F}_\tau(\phi)({\mathcal
L}^{\tau}_\Lambda(f))$ pour tout op\'erateur diff\'erentiel $\phi\in
K[x,d/dx]$. La seconde assertion de la proposition s'ensuit.
\end{proof}

\subsection{5.2 Comparaison de rayons de convergence}\ \\

Soient $\Lambda\in\GL_n(\bQ)$, $\tau$ et $d$ comme au paragraphe
pr\'ec\'edent. Pour toute place $v$ de $K$ et toute matrice $M$ \`a
coefficients dans $K$, on d\'esigne par $\|M\|_v$ le maximum des
valeurs absolues $v$-adiques des \'el\'ements de $M$. On note
d'abord le r\'esultat suivant.\\ \ \\
{\bf Lemme 5.2.1.} \emph{Soit $v$ une place de $K$ au-dessus d'un
nombre premier $p$ qui ne divise pas $d$.  Avec la notation
habituelle $\pi_v = p^{-1/(p-1)}$, on a
\begin{equation}
\label{limites} \lim_{i\to\infty} \|C_{\Lambda,\tau}(i)\|_v^{1/i} =
\pi_v |\tau|_v^{-1}\et \lim_{i\to-\infty}
\|C_{\Lambda,\tau}(i)\|_v^{1/i} = \pi_v^{-1}|\tau|_v.
\end{equation}}
\begin{proof}[Preuve]
Puisque les valeurs propres de $\Lambda$ sont toute rationnelles, il
existe $U\in\GL_n(\bQ)$ telle que le produit $\Delta=U^{-1}\Lambda
U$ soit sous forme de Jordan.  Pour tout $i\in\bZ$, on a
$C_{\Delta,\tau}=U^{-1}C_{\Lambda,\tau}U$ et par suite le rapport
$\|C_{\Lambda,\tau}(i)\|_v/\|C_{\Delta,\tau}(i)\|_v$ est major\'e et
minor\'e par des constantes positives ind\'ependantes de $i$. Cela
permet de supposer que $\Lambda$ est elle m\^eme sous forme de
Jordan.  Dans ce cas, $\Lambda$ est diagonale par blocs avec des
blocs de Jordan $J_1,\dots,J_s$ sur la diagonale et par suite
$C_{\Lambda,\tau}(i)$ est aussi diagonale par blocs avec les blocs
$C_{J_1,\tau}(i), \dots, C_{J_s,\tau}(i)$ sur la diagonale et donc
$\|C_{\Lambda,\tau}(i)\|_v = \max_{1\le k\le s}
\|C_{J_k,\tau}(i)\|_v$. Cela permet de se ramener au cas o\`u
$\Lambda$ est un bloc de Jordan. Dans ce cas, on peut \'ecrire
$\Lambda=r{\bI}_n+N$ o\`u $r\in\bQ\setminus \bZ$ et o\`u $N$ est une
matrice triangulaire sup\'erieure nilpotente. Comme $N^n=0$, on
trouve, pour tout $i\ge 1$,
\begin{equation}
  \label{dec:Cgamma}
\begin{aligned}
C_{\Lambda,\tau}(i)
&= \tau^{-i}\prod_{k=1}^i ((r+k){\bI}_n+N) \\
&= \tau^{-i}(r+1)_i
    \Big( {\bI}_n + \sum_{t=1}^{n-1}
    \Big( \sum_{1\le k_1<\cdots<k_t\le i}
       \frac{1}{(r+k_1)\cdots(r+k_t)}
    \Big) N^t \Big),
\end{aligned}
\end{equation}
Or, pour tout entier $k\ge 1$, la somme $r+k$ est un nombre
rationnel dont la valeur absolue du num\'erateur (au sens usuel) est
major\'ee par $ck$, pour une constante $c>0$ qui ne d\'epend que de
$r$, et dont le d\'enominateur est premier \`a $p$.  On en d\'eduit
$|r+k|_v \ge (ck)^{-1}$ pour tout $k\ge 1$ et par suite, la
d\'ecomposition \eqref{dec:Cgamma} entra\^{\i}ne, pour tout $i\ge
1$,
\begin{equation*}
  |\tau^{-i}(r+1)_i|_v \le \|C_{\Lambda,\tau}(i)\|_v \le (ci)^{n-1} |\tau^{-i}(r+1)_i|_v,
\end{equation*}
la minoration d\'ecoulant du fait que les \'el\'ements de la
diagonale de $C_{\Lambda,\tau}(i)$ sont tous \'egaux \`a $(r+1)_i$.
Cette derni\`ere observation livre aussi $\det C_{\Lambda,\tau}(i) =
\tau^{-in}(r+1)_i^n$. Gr\^ace au lemme 2.2.2, on en d\'eduit
\begin{equation}
\label{est:norme:Cgamma} \lim_{i\to\infty}
\|C_{\Lambda,\tau}(i)\|_v^{1/i} = \lim_{i\to\infty}
|\tau^{-i}(r+1)_i|_v^{1/i} =\pi_v|\tau|_v^{-1},
\end{equation}
ce qui d\'emontre la premi\`ere moiti\'e de \eqref{limites}, et
aussi
\begin{equation}
\label{est:det:Cgamma} \lim_{i\to\infty} \|\det
C_{\Lambda,\tau}(i)\|_v^{1/i} = \lim_{i\to\infty}
|\tau^{-i}(r+1)_i|_v^{n/i} =\pi_v^n|\tau|_v^{-n}.
\end{equation}
On note par ailleurs que, pour toute matrice $A \in\GL_n(K_v)$, on a
\begin{equation*}
\|A\|_v^{-1} \le \|A^{-1}\|_v \le |\det A|_v^{-1} \|A\|_v^{n-1}.
\end{equation*}
En effet, la relation ${\bI}_n=AA^{-1}$ implique $1=\|{\bI}_n\|_v
\le \|A\|_v \|A^{-1}\|_v$ qui livre la premi\`ere des in\'egalit\'es
ci-dessus.  La seconde d\'ecoule simplement de la formule
$A^{-1}=(\det A)^{-1}\Adj(A)$ o\`u $\Adj(A)$ d\'esigne l'adjointe de
$A$.  Si on applique ces formules avec $A=C_{\Lambda,\tau}(i)$, on
d\'eduit de \eqref{est:norme:Cgamma} et de \eqref{est:det:Cgamma}
que
\begin{equation*}
\lim_{i\to\infty} \|C_{\Lambda,\tau}(i)^{-1}\|_v^{1/i} =
\pi_v^{-1}|\tau|_v.
\end{equation*}
Cette formule s'applique aussi \`a $-\Lambda$ au lieu de $\Lambda$
car les valeurs propres de $-\Lambda$ appartiennent aussi \`a
$\bQ\setminus\bZ$.  Alors, gr\^ace \`a \eqref{CC}, on obtient
\begin{equation*}
\lim_{i\to\infty} \|C_{\Lambda,\tau}(-i)\|_v^{1/i} =
\lim_{i\to\infty} \|\tau\Lambda^{-1}
C_{-\Lambda,\tau}(i-1)^{-1}\|_v^{1/i} = \pi_v^{-1}|\tau|_v,
\end{equation*}
ce qui d\'emontre la seconde partie de \eqref{limites}.
\end{proof} \ \\
{\bf Proposition 5.2.2.} \emph{Soient $m$ un entier positif, $v$ une
place finie de $K$ au-dessus d'un nombre premier qui ne divise pas
$d$, et $A(x)$ une matrice de format $m\times n$ \`a coefficients
dans $K_v[[x,1/x]]$. Supposons qu'il existe des nombres r\'eels $r$
et $R$ avec $0 < r \le R$ tels que les coefficients de  $A(x)$
convergent pour tout $x\in\Kbar_v$ avec $r<|x|_v<R$ o\`u $\Kbar_v$
d\'esigne une cl\^oture alg\'ebrique de $K_v$. Alors, on a
${\mathcal L}^{\tau}_\Lambda(A(x)x^\Lambda) = B(x) x^{-\Lambda}$
o\`u $B(x)$ est matrice de format $m\times n$ \`a coefficients dans
$K_v[[x,1/x]]$ dont les coefficients convergent pour tout
$x\in\Kbar_v$ avec $\pi_v |\tau|_v^{-1}R^{-1} < |x|_v <
|\tau|_v\pi_v^{-1} r^{-1}$.}
\begin{proof}[Preuve]
\'Ecrivons $A(x) = \sum_{i=-\infty}^\infty A_i x^i$ de sorte que
$B(x) = \sum_{i=-\infty}^\infty B_i x^{i-1}$ avec $B_i =A_{-i}
C_{\Lambda,{\tau}}(-i)$ pour tout $i\in\bZ$. L'hypoth\`ese que les
coefficients de $A(x)$ convergent en tout point $x$ de $\Kbar_v$
avec $r<|x|_v<R$ signifie que $\limsup_{i\to\infty} \|A_i\|_v^{1/i}
\le 1/R$ et $\limsup_{i\to\infty} \|A_{-i}\|_v^{1/i} \le 1/r$.
Puisque $\|B_i\|_v \le \|A_{-i}\|_v \|C_{\Lambda,\tau}(-i)\|_v$ pour
tout $i\in\bZ$, on en d\'eduit, gr\^ace au lemme 5.2.1, que
$\limsup_{i\to\infty} \|B_i\|_v^{1/i} \le |\tau|_v/(r\pi_v)$ et que
$\limsup_{i\to\infty} \|B_{-i}\|_v^{1/i} \le \pi_v/(R|\tau|_v)$.  La
conclusion suit.
\end{proof}\ \\
\textbf{Remarque.} Cette transform\'ee de Laplace formelle pourrait
avoir d'autres applications, notamment $p$-adiques; Nous comptons
ult\'erieurement l'appliquer aux isocristaux surconvergents.
\section{ Conditions suffisantes}\ \\
Le but de cette section est de compl\'eter la preuve du th\'eor\`eme
1.1. Pour cela, on commence par \'etablir le
r\'esultat pr\'eliminaire suivant:\\ \ \\
\textbf{Proposition 6.1.} \emph{Soit $\psi$ un op\'erateur de
$K[x,d/dx]$ d'ordre $m$, r\'egulier en $0$, ayant une matrice de
r\'eduction $Y_{\psi}\in {\GL}_{m}(K((x)))$ en $0$ telle que
$\prod_{v\in V_{0}}\min(R_{v}(Y_{\psi})\pi_{v}^{-1},1)\ne 0$.
Supposons que $\phi=\overline{\mathcal F}(\psi)$ soit d'ordre $n$ et
r\'egulier en l'infini avec des exposants rationnels non entiers.
Alors $\psi$ est un $E$-op\'erateur.}\\
\ \\
La preuve de cette proposition repose sur les deux lemmes
suivants:\\ \ \\
\textbf{Lemme 6.2.}\emph{ Sous les hypoth\`eses de la proprosition
6.1, l'op\'erateur $\phi$ poss\`ede une base de solutions en
l'infini de forme
$(y_1(\frac{1}{x}),\dots,y_n(\frac{1}{x}))(\frac{1}{x})^{\Gamma}$
telle que $y_1,\dots,y_n\in K[[x]]$, $\Gamma\in{\M}_n({\bQ})$ et
$\prod_{v\in V_{0}}\min(r_{v}(y_i),1)\ne 0$ pour
$i=1,\ldots,n$.}\\ \ \\
\emph{ Preuve.} La d\'ecomposition de Turrittin-Levelt affirme que
les op\'erateurs diff\'erentiels $\phi$ et $\psi$ admettent
respectivement des bases de solutions en l'infini et en $0$ de la
forme
$$
\Big(y_{1}(\frac{1}{x}),y_{2}(\frac{1}{x}),\ldots,y_{n}(\frac{1}{x})\Big)(\frac{1}{x})^{\Gamma}\;\;\;\;\Big(\text{resp.}\;\;\;(w_1,\ldots,w_n):=(F_{1},F_{2},\ldots,F_{m})x^{\Lambda}\Big)$$
telles que \\
1) $y_{1},\ldots,y_{n},F_{1},\ldots,F_{m}\in K[[x]]$,\\
2) $\Gamma=D+N$ est une matrice carr\'ee d'ordre $n$ (resp.
$\Lambda=\Delta+\Theta$ d'ordre
  $m$) o\`u $D$ (resp. $\Delta$) est une matrice diagonale \`a
  coefficients $D_{ii}:=\alpha_{i}$ (resp.
  $\Delta_{ii}:=\beta_{i}$) qui sont les exposants de $\phi$ en
  l'infini (resp. de $\psi$ en $0$),\\
3) $N=(N_{ij})$ (resp.
  $\Theta=(\Theta_{ij})$) est une matrice nilpotente triangulaire
  sup\'erieure \`a coefficients
dans ${\bQ}$ telle que $DN=ND$ (resp.
$\Delta\Theta=\Theta\Delta$).\\

Comme les $\alpha_1,\ldots,\alpha_{n}$ sont des nombres rationnels
non entiers, la proposition 5.1.1 implique que les coefficients du
$n$-uplet
\begin{eqnarray}(\zeta_1,\ldots,\zeta_n)&:=&
{\mathcal
L}_{-\Gamma}^1\Big(\Big(y_{1}(-1/x),y_{2}(-1/x),\ldots,y_{n}(-1/x)\Big)\Big(\frac{1}{x}\Big)^{\Gamma}\Big)
\end{eqnarray}
sont des solutions de $\psi$ en $0$. De plus, le vecteur-ligne
$(\zeta_1,\ldots,\zeta_n)$ peut s'\'ecrire sous la forme
  \begin{eqnarray}
(\zeta_1,\ldots,\zeta_n)&=&(f_1,\ldots,f_n)x^{\Gamma},\;\;\;\text{o\`u}\;\;\;f_1,\ldots,f_n\in
K((x))).
  \end{eqnarray}
Or,
\begin{eqnarray*}
x^{\Gamma}&=& x^{D+N}=\exp(D\log x)\exp(N\log x)=x^{D}\sum_{k\ge
0}\disp\frac{N^{k}}{k!}\log^{k}x \\
           &=&
           x^{D}+x^{D}\sum_{k=1}^{n}\disp\frac{N^{k}}{k!}\log^{k}x.
\end{eqnarray*}
Alors,
\begin{eqnarray}
\zeta_{1}=f_{1}x^{\alpha_{1}},\;\text{et}\;\;\;\;\zeta_{i}&=&f_{i}x^{\alpha_{i}}+\sum_{j=1}^{
i-1}f_{j}x^{\alpha_{j}}\sum_{k=1}^{n}\disp\frac{(N^{k})_{ji}}{k!}\log^{k}x;\;\;\;\;\;2\le
i\le n,
\end{eqnarray}
car $(N^{k})_{ji}=0$ quel que soit $j\ge i$. De m\^eme on trouve
\begin{eqnarray*} w_{1}=F_{1}x^{\beta_{1}}
,\;\;\text{et}\;\;\;\;w_{i}=F_{i}x^{\beta_{i}}+\sum_{j=1}^{
i-1}F_{j}x^{\beta_{j}}\sum_{k=1}^{m}\disp\frac{(\Theta^{k})_{ji}}{k!}\log^{k}x;\;\;\;\;\;2\le
i\le m.
\end{eqnarray*}
Donc,
$$\zeta_{i}\in\; <w_{1},\ldots,w_{m}
>_{K}\;\;\subseteq\;\;
<F_{1}x^{\beta_{1}},\ldots,F_{m}x^{\beta_{m}}>_{K[\log
x]}\;\;(i=1,\ldots,n).$$ En particulier,
$$f_1x^{\alpha_{1}}\;\in\;<F_{1}x^{\beta_{1}},\ldots,F_{m}x^{\beta_{m}}>_{K[\log
x]}.$$ Et par r\'ecurrence sur $i$, on d\'eduit de (6.3) que
$$f_{i}x^{\alpha_{i}}\in
\Big<F_{j}x^{\beta_{j}},f_{k}x^{\alpha_{k}};\;  0\le j\le m, 1\le
k\le i-1\Big>_{K[\log x]}\;\;\;\;\;(i=2,\ldots,n),$$ ou encore,
$$\alpha_{i}\in\{\beta_{j}+h\;|\; h\in{\bZ},\;
j=1,\ldots,m\},\;\;\text{et}\;\;\;\; f_i \in
<F_{1},\ldots,F_{m}>_{K[x,1/x]};\;\;\;\; i=1,\ldots,n.$$ Ce qui
entra\^{i}ne que, pour toute toute place finie $v\in V_0$,
\begin{eqnarray}\min_{1\le i\le n}r_v(f_i)\ge \min_{1\le j\le m}r_v(F_j).
\end{eqnarray} Par ailleurs,
d'apr\`es le lemme 4.4, l'op\'erateur $\psi$ poss\`ede une matrice
de r\'eduction $Y_0$ en $0$ telle que les coefficients de $Y_0^{-1}$
sont des combinaisons $K[x]$-lin\'eaires de $F_1,\ldots,F_m$ et de
leurs d\'eriv\'ees. Ce qui implique, en vertu  du lemme 4.1, que
pour toute place finie $v\in V_0$, $$ \min_{1\le j\le m}r_v(F_j)\ge
r_v(Y_0^{-1})\ge R_v(Y_0^{-1})=R_v(Y_\psi),
$$
et gr\^ace \`a (6.4), on obtient
\begin{eqnarray}
\min_{1\le i\le n}r_v(f_i)\ge R_v(Y_\psi)
\end{eqnarray}
pour toute place finie $v\in V_0$. Or, si on applique la
premi\`ere \'egalit\'e de (5.4) sur\\
$(y_1(\frac{1}{x}),\dots,y_n(\frac{1}{x}))(\frac{1}{x})^{\Gamma}$,
on trouve, depuis (6.1) et (6.2), $${\mathcal
L}_{\Gamma}^1((f_1,\dots,f_n)x^{\Gamma})=-(y_1(-1/x),\dots,y_n(-1/x))\Gamma^{-1}(\frac{1}{x})^{\Gamma}.$$
Ce qui entra\^{i}ne, d'apr\`es la proposition 5.2.2, $\min_{1\le
i\le n} r_v(y_i)\ge \min_{1\le i\le n} r_v(f_i)\pi_{v}^{-1}$ pour
presque tout $v\in V_0$. Par suite, le lemme r\'esulte de
l'hypoth\`ese $\disp\prod_{v\in
V_{0}}\min(R_{v}(Y_\psi)\pi_v^{-1},1)\ne 0.$
\qed \ \\
\ \\
\textbf{Lemme 6.3.}\emph{ Sous les hypoth\`eses de la proprosition
6.1, l'op\'erateur $\phi^*$ poss\`ede lui aussi une base de
solutions en l'infini de forme
$(z_1(\frac{1}{x}),\dots,z_n(\frac{1}{x}))(\frac{1}{x})^{\Lambda}$
telle que $z_1,\dots,z_n\in K[[x]]$, $\Lambda\in{\M}_n({\bQ})$ et
$\prod_{v\in V_{0}}\min(r_{v}(z_i),1)\ne 0$ pour
$i=1,\ldots,n$.} \\ \ \\
\emph{ Preuve.} En vertu du lemme 6.2, il suffit de montrer que les
op\'erateurs $\phi^{*}$ et $\overline{\psi^{*}}$ satisfont les
hypoth\`eses de la proposition 6.1 au lieu  de $\phi$ et $\psi$.
Pour cela, on note d'abord, d'apr\`es \S2.2, que la matrice
$^{t}Y_{\psi}^{-1}$ est une matrice de r\'eduction de
$\mathcal{M}_{\psi^{*}}$ en $0$. Donc, par le changement de variable
$x\longrightarrow -x$, l'op\'erateur $\overline{\psi^{*}}$ est
r\'egulier en $0$ et $Y(x):=^{t}Y_{\psi}^{-1}(-x)$ est une matrice
de r\'eduction de $\mathcal{M}_{\overline{\psi^{*}}}$ en $0$ qui
v\'erifie aussi $\prod_{v\in V_{0}}\min(R_{v}(Y)\pi_{v}^{-1},1)\ne
0$. Donc, gr\^ace au lemme 4.4, si $Y_{\overline{\psi^{*}}}$
d\'esigne une matrice de r\'eduction de $\overline{\psi^{*}}$ en
$0$, alors elle v\'erifie  $\prod_{v\in
V_{0}}\min(R_{v}(Y_{\overline{\psi^{*}}})\pi_{v}^{-1},1)\ne 0$. Par
ailleurs, gr\^ace encore \`a \S2.2, l'op\'erateur $\phi^{*}$ est
r\'egulier en l'infini \`a exposants rationnels non entiers. Enfin,
on a ${\mathcal F}(\phi^{*})=({\mathcal
F}\overline{\phi})^{*}=\overline{\psi^{*}}$ (cf. $[$Ma,
V.3.6$]$).\qed \ \\ \ \\
\emph{Preuve de la proposition 6.1.} En combinant les lemmes 6.2 et
6.3 avec le lemme 4.4 (appliqu\'e en l'infini), on obtient que
$\phi$ poss\`ede une matrice de r\'eduction $Y_\phi$ en l'infini
telle que $ \prod_{v\in V_{0}}\min(R_{v}(Y_\phi),1)\ne 0.$ Comme les
exposants de $\phi$ en l'infini sont rationnels, la proposition
3.2.1 entra\^{i}ne donc que $\phi$ est un $G$-op\'erateur et par
suite $\psi$ est un $E$-op\'erateur.\qed \ \\

Pour compl\'eter la preuve du th\'eor\`eme 1.1, il suffit, d'apr\`es
\S2.5, le th\'eor\`eme 2..4.2 et le corollaire 2.3.3,
de d\'emontrer le th\'eor\`eme suivant: \\
\
\\
{\bf Th\'eor\`eme 6.4.}\emph{ Soit $\psi\in K[x,d/dx]$ un
op\'erateur diff\'erentiel et soit $Y_{\psi}$  une matrice de
r\'eduction de $\psi$ en $0$. Supposons que $\psi$
v\'erifie les conditions suivantes :}\\
(1) \emph{les pentes de $NR(\psi)$ sont dans $\{-1,0\}$;}\\
(2) \emph{les exposants de $\psi$ en $0$ sont rationnels;}\\
(3) $\prod_{v\in V_{0}}\min(R_{v}(Y_{\psi})\pi_{v}^{-1},1)\ne
0$.\\
\emph{Alors $\psi$ est un $E$-op\'erateur.}\\

La principale difficult\'e rencontr\'ee pour prouver ce th\'eor\`eme
vient du fait que l'op\'erateur $\phi:=\overline{\mathcal F}(\psi)$
(r\'egulier en l'infini d'apr\`es corollaire 2.6.2) peut admettre
des exposants entiers en l'infini. Cependant, \`a l'aide des
r\'eductions ci-dessous, on se ram\`ene au cas o\`u tous les
exposants de $\phi$ en l'infini sont rationnels non entiers et alors
le th\'eor\`eme
d\'ecoule de la proposition 6.1.\\ \ \\
{\bf R\'eductions}\\
I. Tout d'abord, on remarque que, pour tout $\ell\in \mathbb{Z}$ tel
que $x^{\ell}\psi\in K[x,d/dx]$, l'op\'erateur $x^{\ell}\psi$
v\'erifie les conditions du th\'eor\`eme. En effet, le polygone de
Newton Ramis de $x^{\ell}\psi$ n'est autre que celui de $\psi$ \`a
une translation verticale pr\`es. Les conditions (2) et (3) se
d\'eduisent du fait que $A_{\psi}=A_{x^{\ell}\psi}$ et  ${\mathcal
M}_{\psi}={\mathcal M}_{x^{\ell}\psi}$. Donc, d'apr\`es (1),
corollaire 2.6.2 et lemme 3.3.1, il suffit de d\'emontrer le
th\'eor\`eme pour les op\'erateurs de la forme:
$\psi=\sum_{i=0}^{m}a_{i}(x)\Big(x\disp\frac{d}{dx}\Big)^{i}$, avec
$ a_{i}(x)=\sum_{j=0}^{n}a_{ij}x^{j}\in K[x]$ et $
a_{m}(x)=1.$\\
II. Maintenant, si $\phi$ a des exposants entiers en l'infini, on
consid\`ere les op\'erateurs diff\'erentiels suivants d\'eduits de
$\psi$ en fonction d'un param\`etre $\alpha\in\mathbb{Q}$:
$$\begin{array}{ccccccc}
    &   &   & \overline{\mathcal F} &   &   &   \\
  \phi &=&
\disp\sum_{i=0}^{m}a_{i}\Big(\disp\frac{d}{dx}\Big)\Big(-x\disp\frac{d}{dx}-1\Big)^{i}
                  & \longleftarrow     &
   \psi                     & = &
\disp\sum_{i=0}^{m}a_{i}(x)\Big(x\disp\frac{d}{dx}\Big)^{i} \\
       & &  & &
        \downarrow & &\\
_{\alpha}\phi &=
&\disp\sum_{i=0}^{m}a_{i}\Big(\disp\frac{d}{dx}\Big)\Big(-x\disp\frac{d}{dx}-\alpha-1\Big)^{i}
   & \longleftarrow     &
  \psi_{\alpha}            & =&
\disp\sum_{i=0}^{m}a_{i}(x)\Big(x\disp\frac{d}{dx}-\alpha\Big)^{i}.\\
&   &   & \overline{{\mathcal F}} &   &   &
\end{array}
$$
L'op\'erateur $\psi_{\alpha}$ ainsi d\'efini  partage avec $\psi$
les propri\'et\'es suivantes.
\begin{itemize}
\item[(1)] Les pentes de  $NR(\psi_{\alpha})$ appartiennent \`a
$\{-1,0\}$. \item[(2)] Tous les exposants de $\psi_{\alpha}$ en $0$
sont rationnels. \item[(3)] $\prod_{v\in
V_{0}}\min(R_{v}(Y_{\psi_{\alpha}})\pi_{v}^{-1},1)\ne 0$, o\`u
$Y_{\psi_{\alpha}}$ est une matrice de r\'eduction de
$\psi_{\alpha}$ en $0$.
\end{itemize}
En effet, la premi\`ere assertion r\'esulte du lemme 2.5.1. Quant
aux deux derni\`eres, on va les d\'emontrer simultan\'ement. La
matrice $$ A= \frac{1}{x}\left (
\begin{array}{c}
0\\0\\\vdots\\0\\-a_0
\end{array}
\begin{array}{c}
1\\0\\\vdots\\0\\-a_1
\end{array}
\begin{array}{c}
0\\1\\\vdots\\0\\-a_2
\end{array}
\begin{array}{c}
\ldots\\ \ldots\\ \ddots \\ \ \\ \ldots
\end{array}
\begin{array}{c}
0\\0\\\vdots\\1\\-a_{m-1}
\end{array}
\right )
$$ est associ\'ee  au $K(x)$-module diff\'erentiel $\mathcal{M}_\psi$.
En combinant les hypoth\`eses (1) et (2) du th\'eor\`eme 6.4 avec
les observations de \S2.5, la matrice $A$ poss\`ede une matrice de
r\'eduction $Y_A\in {\GL}_m(K((x)))$ en $0$  telle que
$Y_A[A]=x^{-1}\Lambda$ o\`u  $\Lambda$ d\'esigne une matrice
$m\times m$ \`a coefficients dans $K$ dont les valeurs propres sont
dans ${\bQ}$. Donc
$Y_A[A+\disp\frac{\alpha}{x}{\bI}_m]=\disp\frac{1}{x}(\Lambda+\alpha
{\bI}_m)$. Par cons\'equent, la matrice $Y_A$ est aussi une matrice
de r\'eduction de $\Big(A+\disp\frac{\alpha}{x}{\bI}_m\Big)$. Or
$\Big(A+\disp\frac{\alpha}{x}{\bI}_m\Big)$ est une matrice
associ\'ee au $K(x)$-module diff\'erentiel ${\mathcal
M}_{\psi}\otimes_{K(x)}{\mathcal M}_{\alpha}= {\mathcal
M}_{\psi_{\alpha}} $. Donc, d'une part, les exposants de
$\psi_{\alpha}$ sont tous rationnels en $0$. D'autre part, depuis
(3) du th\'eor\`eme 6.4, le corollaire 4.3, appliqu\'e \`a $Y_\psi$
et $Y_A$ ensuite \`a $Y_A$ et $Y_{\psi_{\alpha}}$, entra\^{i}ne
$\prod_{v\in V_{0}}\min(R_{v}(Y_{\psi_{\alpha}})\pi_{v}^{-1},1)\ne
0$, et la
conclusion suit.\\

Par ailleurs, les formules suivantes
\begin{eqnarray*}
_{\alpha}\phi
(x^{-s})&=&\sum_{i=0}^{m}a_{i}\Big(-\disp\frac{d}{dx}\Big)\Big(x\disp\frac{d}{dx}-\alpha-1\Big)^{i}(x^{-s})\\
                  &=&\sum_{i=0}^{m}a_{i0}x^{-s}(s-\alpha-1)^{i}+
                  \mbox{(termes d'ordre inf\'erieur)},\\
\end{eqnarray*}
et
\begin{eqnarray*}
\psi_{\alpha}
(x^{s})&=&\sum_{i=0}^{m}a_{i}(x)\Big(x\disp\frac{d}{dx}-\alpha\Big)^{i}(x^{s})\\
                  &=&\sum_{i=0}^{m}a_{i0}x^{s}(s-\alpha)^{i}+
                  \mbox{(termes d'ordre sup\'erieur)},\\
\end{eqnarray*}
montrent que les polynomes indiciels  des op\'erateurs de
$_{\alpha}\phi$ et $\phi$ en l'infini sont respectivement
$$\sum_{i=0}^{m}a_{i0}(s-\alpha-1)^{i}\;\;\;\;\;\;\mbox{et}\;\;\;\;\;
\sum_{i=0}^{m}a_{i0}(s-1)^{i},$$  et que ceux de $\psi_{\alpha}$ et
$\psi$ en $0$ sont respectivement
$$\sum_{i=0}^{m}a_{i0}(s-\alpha)^{i}\;\;\;\;\;\;\mbox{et}\;\;\;\;\;
\sum_{i=0}^{m}a_{i0}s^{i}.$$ Donc, d'une part, les exposants de
$_{\alpha}\phi$ en l'infini sont ceux de $\phi$ en l'infini
translat\'es par  $-\alpha$, et  d'autre part, les exposants de
$_{\alpha}\phi$ en l'infini sont ceux de $\psi_{\alpha}$ en $0$
translat\'es par  $1$. Cela montre que les exposants de
$_{\alpha}\phi$  et $\phi$ en l'infini sont tous rationnels. De
plus, comme  les op\'erateurs $\phi$ et $_{\alpha}\phi$ sont
r\'eguliers en l'infini, ils ont  le m\^eme ordre $m$ que celui de
$\psi$ (cf. $[$Ba, \S2. Remark 6$]$). Par cons\'equent, quitte \`a
choisir un rationnel convenable $\alpha$, les exposants de
$_{\alpha}\phi$ en l'infini sont tous des nombres rationnels non
entiers et le th\'eor\`eme r\'esulte de la proposition 6.1.\qed

\end{document}